\documentclass[a4paper,12pt]{article}
\pdfoutput=1
\usepackage[utf8]{inputenc}
\usepackage{amsfonts}
\usepackage{graphicx,color,epstopdf}
\usepackage{amsmath,amsthm, morefloats, bm}
\allowdisplaybreaks

\def\Xint#1{\mathchoice
{\XXint\displaystyle\textstyle{#1}}%
{\XXint\textstyle\scriptstyle{#1}}%
{\XXint\scriptstyle\scriptscriptstyle{#1}}%
{\XXint\scriptscriptstyle\scriptscriptstyle{#1}}%
\!\int}
\def\XXint#1#2#3{{\setbox0=\hbox{$#1{#2#3}{\int}$}
\vcenter{\hbox{$#2#3$}}\kern-.5\wd0}}
\def\cint{\Xint-}

\def\d{{\rm dist}}
\begin{document}
\title{Perfectly-matched-layer boundary integral equation
method for wave scattering in a layered medium}
\author{{Wangtao Lu}
	\thanks{Department of Mathematics, Michigan State University, East
	Lansing, MI 48824, USA.  Email: wangtaol@math.msu.edu}
	\and{Ya Yan Lu}
	\thanks{Department of Mathematics, City University of Hong Kong, Kowloon,
	Hong Kong. Email:mayylu@cityu.edu.hk}
	\and{Jianliang Qian}
	\thanks{Department of Mathematics, Michigan State University, East
	Lansing, MI 48824, USA.  Email: qian@math.msu.edu}
}
\maketitle
\begin{abstract}
For scattering problems of time-harmonic waves, the boundary integral
equation (BIE) methods are highly competitive, since they are formulated on
lower-dimension boundaries or interfaces, and can automatically satisfy
outgoing radiation conditions. For scattering problems in a layered
medium, standard BIE methods based on the Green's function of the
background medium must evaluate the expensive Sommefeld integrals.
Alternative BIE methods based on the free-space Green's function give rise
to integral equations on unbounded interfaces which are not easy to
truncate, since the wave fields on these interfaces decay very slowly. We
develop a BIE method based on the perfectly matched layer (PML) technique.
The PMLs are widely used to suppress outgoing waves in numerical methods
that directly discretize the physical space. Our PML-based BIE method uses
the Green's function of the PML-transformed free space to define the
boundary integral operators. The method is efficient, since the Green's
function of the PML-transformed free space is easy to evaluate and the
PMLs are very effective in truncating the unbounded interfaces. Numerical
examples are presented to validate our method and demonstrate its accuracy.
\end{abstract}
%\tableofcontents
%\newpage
\section{Introduction}
Scattering problems for sound, electromagnetic and elastic waves in
layered media are highly relevant for practical applications \cite{chew95}. Numerical methods that directly discretize the physical domain,
such as the finite element method (FEM) \cite{mon03}, are very
versatile and widely used, but they become too expensive when the scatterer
is large compared with the wavelength. The boundary integral equation
(BIE) methods \cite{colkre13} are applicable to structures with
piecewise constant material parameters. These methods take care of the
outgoing radiation condition automatically and reduce the dimension by
one, since the integral equations are formulated on material
interfaces or boundaries of obstacles. For many problems, BIE methods
can outperform FEM and other domain-discretization methods, and
deliver highly accurate solutions with relatively small computing efforts. 

For scattering problems in a layered medium, the common
BIE methods are based on the Green's function of the layered
background medium \cite{derpol35,som1909,weyl1919}, 
so that the intergral equations are formulated on strictly local
interfaces or boundaries. However, it is well known
that this approach is bottlenecked by the evaluation of Sommefeld
integrals arising from the layered-medium Green's function and its
derivatives. Over the past decades, many methods
such as high-frequency asymptotics,
rational approximations, contour deformations 
\cite{cai13,caiyu00,okhcan04,paugaybalmar00,perarabru14}, complex images
\cite{och04,tar05,thowea75}, and the steepest descent method
\cite{cuiche98,cuiche99}, have been developed to speed up the
computation of Sommefeld integrals.  
Unfortunately, the computational cost for
evaluating the Sommerfeld integrals remains high \cite{cai02}.

An alternative approach is to use the free-space Green's function, but
then the integral equations must also be formulated on the
unbounded interfaces separating the different layers of the background
medium.  Various types of compactly supported functions can be used to truncate the 
unbounded interfaces and to suppress the artifical reflections from 
the edges of the truncated sections. 
Existing methods in this category include the approximate truncation method
\cite{meicha01,saisor11}, the taper function method
\cite{zhalismicar05,spisorsai08,mirsorsai14}, and the windowing 
function method \cite{brudel14,mon08,brulyoperaratur16,laigreone15}.
In particular,  the windowing function method of Bruno {\it et al.}
\cite{brulyoperaratur16} can largely eliminate the artificial
reflections, since the errors decrease superalgebraically as the
window size is increased. 
Similar good performance can be observed in the hybrid method of Lai
{\it et al.} \cite{laigreone15} that combines windowed layer
potentials (in physical space) with a Sommerfeld-type correction (in
Fourier space) for scattering problems where the obstacles are close
to or even cut through the interfaces of the background layered media. 
 
In this paper, we develop a BIE method based on perfectly matched
layers (PMLs) for two-dimensional (2D) scattering problems in layered
media. The PML technique is widely used for domain truncations
in wave propagation problems \cite{ber94, ber96, chexia13, lassom01}. It can
be regarded as a complex coordinate stretching that replaces the real
independent variables in the original governing equation by complex
independent variables, so that the outgoing waves are damped as they
propagate into the PML region. Similar to those BIE methods based on the
free-space Green's function,  our BIE method avoids evaluating the expensive
Sommefeld integrals, but requires integral equations along the interfaces of
the background layered medium. But instead of the free-space Green's
function, we use the Green's function for the PML-transformed free space, so
that the truncation of the interfaces follows automatically from the
truncation of PMLs. Notice that the Green's function of the PML-transformed
free space can be simply obtained by extending the argument of the usualy
Green's function to complex space following the definition of the complex
square root function. 

We implement our PML-based BIE method for 2D scattering problems involving two
homogeneous media separated by a single interface. The interface is flat
except in a finite session which is referred to as the local
perturbation. Additional obstacles are also allowed in the homogeneous
media. Two common types of incident waves are considered: a plane incident wave
and  a cylindrical wave due to a point source. The integral
equations are established for the scattered wave 
satisfying Sommefeld radiation condition at infinity. 
The scattered wave is defined as the difference between the total wave field
and a reference wave field obtained from the same indicent wave for
the layered background medium (without the local perturbation of the
interface and the obstacles). 

BIE methods for scattering problem use many different formulations.
Some of these formulations are more appropriate for large (i.e.
high-frequency) problems, since they give rise to linear systems with
better condition numbers, and are thus more efficient when iterative
methods are used. Since our purpose is to demonstrate the
effectiveness of PML-based BIEs for truncating unbounded interfaces,
we adopt a simple formulation that comes from Green's representation
theorem directly. In addition, we calculate the so-called
Neumann-to-Dirichlet (NtD) map (mapping Neumann data to Dirichlet data
on the boundary) for each subdomain with constant material parameters,
so that the final linear system on interfaces or boundaries of the
obstacles can be written down in a very simple form.

To approximate the integral equations, we utilize a graded
mesh technique \cite{colkre13}, a high-order quadrature rule by Alpert
\cite{alp99}, and a newly proposed stabilizing technique. 
Numerical results indicate that our method is highly accurate and the
truncation of the unbounded interfaces by PML is very effective. 
Typically, for a PML with a thickness of one wavelength and
discretized with about the same number of points as a typical segment
of one wavelength, about seven significant digits can be obtained. Numerical
results show that numerical error decays exponentially for $S$ (a PML
parameter representing the strength of the PML) in whatever range. 

The rest of this paper is organized as follows. In sections 2 and 3, we
present our PML-based BIE formulation for solving scattering problems in
layered media. 
The numerical schemes for discretizing the integral equations are
given in sections 4 and 5. 
Numerical examples are presented in section 6 to illustrate
the performance of our method, and we conclude the paper in section 7.

\section{Problem Formulation}
In this paper, we mainly focus on two-dimensional TE and TM polarized
scattering problems in a planar layered medium with local perturbations
and/or obstacles. To clarify our methodology in a simpler setting, we assume
that only local perturbations exist in the medium in the following. In
general, obstacles make no noticeable difficulties for the scattering
problem.  

\begin{figure}[!ht]
  \centering
    \includegraphics[width=0.6\textwidth]{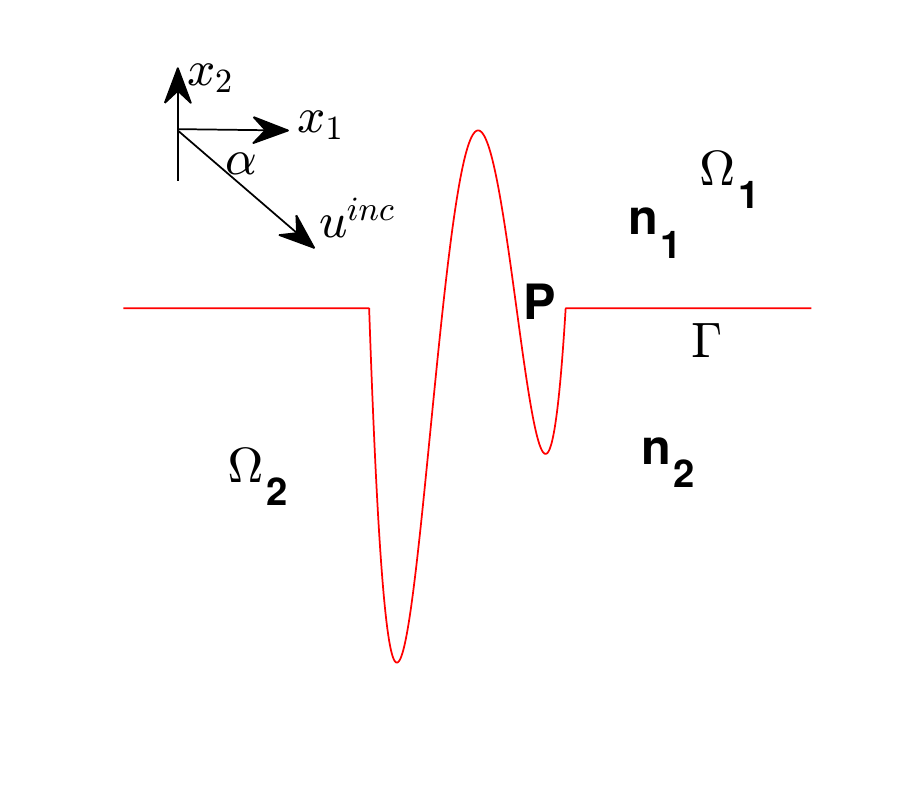}
  \caption{Profile of a 2D layered medium.}
		\label{fig:profile1}
\end{figure}

As illustrated in Figure~\ref{fig:profile1},
the layered medium is $x_3$-invariant and consists of two homogeneous layers
$\Omega_j$ with constant refractive index $n_j$ for $j=1,2$. The interface
$\Gamma$ separating the two layers is flat on $x_2=0$ but contains a local
perturbation curve $P$, smooth or piecewise smooth. Here, $(x_1,x_2,x_3)$
denotes the standard Cartesian coordinate system.   

The total field $u^{tot}$, representing the $x_3$-component of
electric field in TE polarization or the $x_3$-component of magnetic field
in TM polarization, solves 
\begin{align}
	\label{eq:gov:u}
	\Delta u^{tot} + k_0^2n_j^2 u^{tot} &= 0,\quad{\rm in}\quad\Omega_j,\\
	\label{eq:trans:cond}
	[u^{tot}]=0, \quad\left[\frac{\eta_j\partial u^{tot}}{\partial {\bm \nu}}\right] &=
	0,\quad{\rm on}\quad\Gamma,  
\end{align}
where $k_0=\frac{2\pi}{\lambda}$ is the freespace wavenumber, $\lambda$
is the wavelength, ${\bm \nu}$ denotes the unit normal vector along $\Gamma$
pointing toward $\Omega_2$, $[f]$ denotes the jump of the quantity $f$
across $\Gamma$, $\eta_j=1$ in TE polarization and $\eta_j=\frac{1}{n_j^2}$
in TM polarization. Let $u^{inc}$ be an incident wave from the upper medium
$\Omega_1$, and then one usually rewrites
\begin{equation}
	\label{eq:define:ut}
	u^{tot} = \left\{
\begin{array}{lc}
	u^{inc} + u^{r}_1, \quad{\rm in}\quad\Omega_1,\\
	u^{t}_2, \quad{\rm in}\quad\Omega_2,
\end{array}
		\right.  
\end{equation}
where $u^{r}_1$ represents the reflective wave in $\Omega_1$ and $u^{t}_2$
represents the transmitted wave in $\Omega_2$.  

In the following, we focus on two common types of incident waves, a
plane incident wave and a cylindrical wave due to a source
$x^*=(x_1^*,x_2^*)\in\Omega_1$. In the latter case, equation
(\ref{eq:gov:u}) should be replaced by 
\begin{equation}
	\label{eq:gov:u:r}
	\Delta u^{tot} + k_0^2n_j^2 u^{tot} = -\delta(x,x^*),\quad{\rm
	in}\quad\Omega_j,
\end{equation}
so that the total field $u^{tot}$ represents a layered-medium Green's
function at the source $x^*$.

We first discuss the case for plane incident waves. Suppose the incident wave
is given by $u^{inc} = e^{ik_0n_1(x_1 \cos\alpha - x_2 \sin\alpha )}$, where
$\alpha\in[0,\pi]$ denotes the angle between the wave direction and the
positive $x_1$-axis. Neither $u^r_1$ nor $u^t_2$ satisfies the
Sommerfeld radiation condition since neither of them is outgoing
in all directions. To extract an outgoing wave field, we need a reference
solution, denoted by $u^{tot}_0$, to the scattering problem with perfectly
flat interface $x_2=0$ and with the same incident wave $u^{inc}$. One easily
gets that
\begin{equation}
	\label{eq:u0}
u_0^{tot} = \left\{
\begin{array}{lc}
	e^{ik_0n_1(x_1 \cos\alpha - x_2 \sin\alpha )} + (T-1)e^{ik_0n_1(x_1\cos\alpha 
	 +x_2\sin\alpha )}, & {\rm in}\quad\Omega_1,\\
		 T e^{ik_0n_1x_1\cos\alpha  - i k^* x_2}, & {\rm in}\quad\Omega_2,
\end{array}
\right.  
\end{equation}
where 
\begin{align*}
	k^* &= k_0\sqrt{n_2^2-n_1^2\cos^2\alpha},\\
	T &= \frac{2}{1+\frac{k^*\eta}{k_0n_1\sin\alpha}},
\end{align*}
and $\eta=\eta_1/\eta_2$.  
Then, 
\begin{equation}
	\label{eq:def:us}
	u^s = u^{tot}-u_0^{tot}:=\left\{
		\begin{array}{lc}
			u^s_{1} & {\rm in}\quad\Omega_1,\\
			u^s_{2} & {\rm in}\quad\Omega_2,
		\end{array}
		\right.
\end{equation}
defines an outgoing wave that satisfies 
\begin{align}
	\label{eq:gov:us}
	\Delta u^s_j + k_0^2n_j^2 u^s_j &= 0,\quad{\rm in}\quad\Omega_j,\\
	        \label{eq:som:cond}
        \lim_{r\rightarrow\infty}r^{-1/2}\left(\frac{\partial
        u^s_j}{\partial r} - ik_0 n_j u^s_j\right) &=0,\quad{\rm
          in}\quad\Omega_j.  
\end{align}
The transmission condition (\ref{eq:trans:cond}) then becomes
\begin{align}
\label{eq:cond:us1}
	\left.u^s_1\right|_{\Gamma} - \left.u^s_2\right|_{\Gamma} &=-[u_0^{tot}],\\
	\label{eq:cond:us2}
	\left.\eta_1\frac{\partial u^s_1}{\partial {\bm \nu}}\right|_{\Gamma} -
        \left.\eta_2\frac{\partial u^s_2}{\partial {\bm \nu}}\right|_{\Gamma}&=
        -\left[\eta_j\frac{\partial u_0^{tot}}{\partial {\bm \nu}}\right].
\end{align}
We note that away from the local perturbation curve $P$, $u^s = u^s_1=u^s_2$ and
$\eta_1\partial_{{\bm \nu}}u^s_1 = \eta_2\partial_{{\bm \nu}}u^s_2$ on $\Gamma$.   

On the other hand, if the incident wave is a cylindrical wave
$u^{inc} = \frac{i}{4} H_0^{(1)}(k_0n_1|x-x^*|)$ due to a source point
$x^*\in\Omega_1$. In this case, one easily obtains that, by defining
\begin{equation}
	\label{eq:utot0:2}
u_0^{tot} = \left\{
\begin{array}{lc}
	u^{inc}, & {\rm in}\quad\Omega_1,\\
		 0, & {\rm in}\quad\Omega_2,
\end{array}
\right.
\end{equation}
the difference wave field $u^s = u^{tot}-u^{tot}_0$ defines an outgoing wave. 

In a typical BIE formulation, the computation of $u^s$ in the whole plane
can be reduced to computing $u^s_j$ and $\partial_{\bm \nu} u_j^s$ on
$\Gamma$, governed by the transmission conditions (\ref{eq:cond:us1}) and
(\ref{eq:cond:us2}). To solve (\ref{eq:cond:us1}) and
(\ref{eq:cond:us2}), we require a relation between $u_j^s$ and
$\partial_{\bm\nu} u_j^s$ for $j=1,2$. In this paper, we make use of
Neumann-to-Dirichlet maps ${\cal N}_j$ that satisfies $u^s_j = {\cal
N}_j\partial_{{\bm \nu}} u^s_j$ on the boundary $\Gamma$ for each outgoing
wave $u^s_j, j=1,2$. Then, (\ref{eq:cond:us1}) and (\ref{eq:cond:us2})
become
\begin{equation}
  \label{eq:linsys:us}
\left[
\begin{array}{cc}
	{\cal N}_s^1 & - {\cal N}_s^2 \\
	\eta_1{\cal I} & -\eta_2{\cal I} 
\end{array}
	\right]\left[
		\begin{array}{c}
			\left.\partial_{{\bm \nu}} u^s_1\right|_{\Gamma}\\
			\left.\partial_{{\bm \nu}} u^s_2\right|_{\Gamma}
		\end{array}
	\right]
	= \left[
		\begin{array}{c}
		-[{u}_0^{tot}]\\
-\left[\eta_j\partial_{\bm \nu} {u}_0^{tot}\right]
		\end{array}
	\right],
\end{equation}
where ${\cal I}$ denotes the identity operator. On solving
(\ref{eq:linsys:us}), we immediately obtain that $u^s_j|_{\Gamma} = {\cal
N}_s^j \partial_{{\bm \nu}} u^s_j|_{\Gamma}$.  	

In practice, since $\Gamma$ is unbounded and since $u^s_j$ decays slowly at
infinity, it is impossible to find a finite-dimensional matrix to
accurately approximate ${\cal N}_s^j$ by directly discretizing the whole
boundary of $\Gamma$ without truncating it.  To resolve this issue, we use
a PML to enclose the local perturbation curve $P$
so that any outgoing wave can be absorbed. Therefore, local
transmission condition can be imposed on a finite session of $\Gamma$
including $S$.  In doing so, we first need to construct NtD maps for domains
in a PML environment, and this relies on boundary integral equations.    

\section{Boundary integral equation in a half space}
Without loss of generality, we consider the outgoing solution $u^s_1$ in
$\Omega_1$, that satisfies
\begin{align}
	\label{eq:hemlsct}
	\Delta u^s_1 + k_0^2 n_1^2 u^s_1 &= 0, \\
	\label{eq:som:sct}
	\lim_{r\rightarrow \infty}r^{1/2} \left(\frac{\partial u^s_1}{\partial r} -
	ik_0 n_1u^s_1\right) &= 0,\quad r=|{x}|, 
\end{align}
for $x\in\Omega_1$.  To simplify the presentation in this section, we assume
that the piecewise smooth curve $S$ is bounded by a box $[-a_1, a_1]\times
[-a_2, a_2]$ for $a_j>0$. Unless otherwise specified, we will suppress
the subscript $1$ indexing the domain $\Omega_1$ so that we use $\Omega$,
$u^s$, and $n$ to denote $\Omega_1$, $u^s_1$, and $n_1$, respectively.

\subsection{BIE in physical domain}
As is well-known, the fundamental solution to equation (\ref{eq:hemlsct}) is
\begin{equation}
	\label{eq:fund}
	G(x, y) = \frac{i}{4}H_0^{(1)}(k_0n|x-y|),
\end{equation}
which solves
\begin{equation}
	\label{eq:ps:helm}
	\Delta_{x} G(x, y) + k_0^2n^2 G(x,y) = -\delta(x-y),
\end{equation}
for $x,y\in\Omega$.  

\begin{figure}[!ht]
  \centering
    \includegraphics[width=0.7\textwidth]{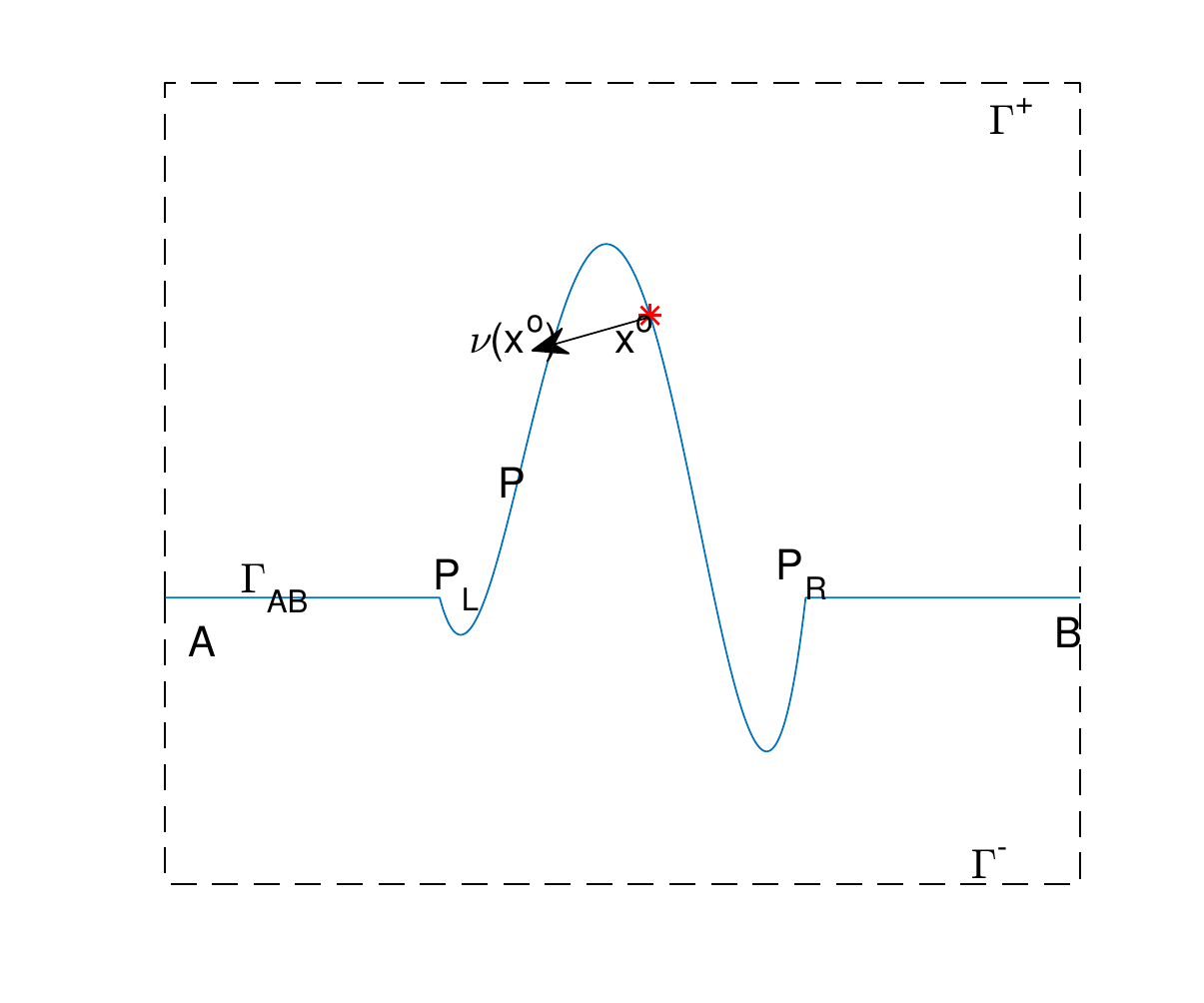}
  \caption{Profile of a 2D layered medium with direct truncation.}
		\label{fig:profile:t}
\end{figure}

As shown in Figure~\ref{fig:profile:t},
to truncate the unbounded interface $\Gamma$, we place a box bounded by 
$\Gamma^+\cup\Gamma^-$ to enclose the local perturbation curve $P$ so that
$\Gamma$ truncated by the box becomes a bounded curve $\Gamma_{AB}$, which
is composed of $AP_L$, $P$, and $P_RB$. Clearly, $\Omega$ is truncated into
a domain $\Omega^s$ bounded by $\Gamma^s=\Gamma_{AB}\cup\Gamma^+$, where
$\Gamma^+$ is the dashed line above $\Gamma_{AB}$.    

According to \cite{colkre13}, one easily obtains the following
representation theorem
\begin{equation}
	\label{eq:rep}
	u^s(x) = \int_{\Gamma^s} \partial_{{\bm \nu}} G(x, y) u^s(y) -
        G(x,y)\partial_{{\bm \nu}} u^s(y)ds(y),   
\end{equation}
for $x\in\Omega^s$. As $x$ approaches $\Gamma^s$, one gets the
following boundary integral equation (see \cite{colkre13,lulu14}) 
\begin{equation}
	\label{eq:bie:org}
	({\cal K}-{\cal K}_0[1]) [u^s](x)= {\cal S} [\partial_{{\bm \nu}} u^s](x),  
\end{equation}
for $x\in\Gamma^s$.  Here, we have defined the following boundary integral
operators
\begin{align}
  \label{eq:def:oS}
	{\cal S}[\phi](x) &= 2\int_{{\Gamma^s}} {G}(x,y) \phi(y) ds(y),\\
  \label{eq:def:oK}
	{\cal K}[\phi](x) &= 2\cint_{{\Gamma^s}} \partial_{{\bm \nu}} {G}(x,y)
  \phi(y) ds(y),\\
  \label{eq:def:oK0}
	{\cal K}_0[\phi](x) &= 2\cint_{{\Gamma^s}} \partial_{{\bm \nu}} {G}_0(x,y)
  \phi(y) ds(y),
\end{align}
where $G_0(x,y)=\frac{1}{2\pi}\log|x-y|$ is the Green's function of
Laplacian equation $\Delta u = 0$, and $\cint$ denotes the Cauchy principal  
integral.  Therefore, one obtains the NtD operator ${\cal N} = ({\cal
K}-{\cal K}_0[1])^{-1}{\cal S}$ that maps $\partial_{{\bm \nu}} u^s$ to $u^s$ on the
bounded curve $\Gamma^s$. 

Now, a significant question arises: what boundary conditions should we
impose on $\Gamma^+$?  One may directly specify that $u^s\approx 0$ and
$\partial_{\bm \nu} u^s\approx 0$ on $\Gamma^+$ to truncate the NtD operator
${\cal N}$ onto $\Gamma_{AB}$.  
Unfortunately, the outgoing wave $u^s(x)$ can decay slowly as $x$ approaches
infinity in $\Omega$. Of course, we may place $\Gamma^+$ sufficiently far
away from the perturbation curve $P$.  However, the computational domain can
become extremely large whereas the boundary condition still maintains a
low-order accuracy.  To address this issue, we propose to introduce a PML to
surround the local perturbation $P$, as will be presented below.  

\subsection{Green's representation theorem in PML-transformed domain}
Specifically, we introduce the complex coordinate stretching function
$\tilde{x}(x)=(\tilde{x}_1(x_1), \tilde{x}_2(x_2))$ by defining 
\begin{align}
	\label{eq:x}
	\tilde{x}_l(x_l) = x_l + i\int_{0}^{x_l} \sigma_l(t) dt,
\end{align}
for $l=1, 2$, where we take 
\begin{align}
	\label{eq:sigma}
	\sigma_l(t)=\sigma_l(-t), \sigma_l=0\,\,{\rm for}\,\, |t|\leq a_l,\,{\rm
	and}\, \sigma_l(t)\geq 0\,\,{\rm for}\,\, |t|\geq a_l.  
\end{align}
Domains with nonzero $\sigma_l$ are called the {\it perfectly matched
layer}. Since $\sigma_l$ is $0$ in $[-a_1,a_1]\times[-a_2,a_2]$, 
the PML does not overlap the local perturbation $P$; the setup of
$\sigma_l(t)$ will be discussed later.  

Based on (\ref{eq:rep}), we can analytically continue $u^s$ onto the domain
$\tilde{\Omega}^s=\{\tilde{x}(x)|x\in\Omega^s\}$ by defining 
\begin{align}
	\label{eq:rep:ext}
	u^s(\tilde{x}) = \int_{\Gamma^s} \partial_{{\bm \nu}} G(\tilde{x}, y) u^s(y) -
	G(\tilde{x},y)\partial_{{\bm \nu}} u^s(y)dy.   
\end{align}

According to \cite{chexia13}, one sees that $u^s(\tilde{x})$ satisfies
\begin{equation}
	\label{eq:pml:helm}
	\tilde{\Delta} u^s(\tilde{x}) + k_0^2n^2 u^s(\tilde{x}) = 0,\quad{\rm
	in}\quad\tilde{\Omega}^s,  
\end{equation}
where $\tilde{\Delta}=\partial_{\tilde{x}_1}^2 + \partial_{\tilde{x}_2}^2$.  
%Outside the PML, $u^s(\tilde{x})=u^s(x)$ since $\tilde{x}_1=x_1$ and
%$\tilde{x}_2=x_2$. Inside the PML, $G(\tilde{x},\cdot)$ and $\partial_{\bm
%\nu} G(\tilde{x}, \cdot)$ exponentially decay to $0$ (see \cite{lassom01})
%as $|x|\rightarrow\infty$ so that $u^s(\tilde{x})$ exponentially decay to
%$0$ as well.   
Defining the complexified function $\tilde{u}^s(x)=u^s(\tilde{x})$ on
$\Omega^s$, we see that equation (\ref{eq:pml:helm}) can be rewritten by the
chain rule as
\begin{equation}
	\label{eq:pml:helm:2}
        \nabla\cdot({\bf A} \nabla \tilde{u}) + k_0^2 n^2J\tilde{u} = 0,
\end{equation}
where $\alpha_1(x_1) = 1+i\sigma_1(x_1)$, $\alpha_2(x_2) =
1+i\sigma_2(x_2)$, ${\bf A}={\rm diag}\{\alpha_2/\alpha_1,
\alpha_1/\alpha_2\}$, and $J(x)=\alpha_1(x_1)\alpha_2(x_2)$.  

As shown in \cite{lassom01}, the fundamental solution to
(\ref{eq:pml:helm:2}) is 
\begin{equation}
	\label{eq:fun:helm2}
	\tilde{G}(x,y) = G(\tilde{x},\tilde{y}) =
	\frac{i}{4} H_0^{(1)}(k_0n\rho(\tilde{x},\tilde{y})),
\end{equation}
where the complexified distance function
$\rho$ is defined to be 
\begin{equation}
	\label{eq:comp:dist}
	\rho(\tilde{x},\tilde{y}) = [(\tilde{x}_1-\tilde{y}_1)^2
	+ (\tilde{x}_2-\tilde{y}_2)]^{1/2},
\end{equation}
and the half-power operator $z^{1/2}$ is chosen to be the branch of
$\sqrt{z}$ with nonnegative real part for $z\in \mathbb{C}/(-\infty,0]$; in
other words, $\tilde{G}$ satisfies  
\begin{equation}
	\label{eq:comp:fund}
	\nabla_y\cdot({\bf A}(y)\nabla_y\tilde{G}(x,y))
	+k_0^2 n^2 J(y)\tilde{G}(x,y)=-\delta_x(y),
\end{equation}
for $x,y\in\Omega^s$.  

\begin{figure}[!ht]
  \centering
    \includegraphics[width=0.7\textwidth]{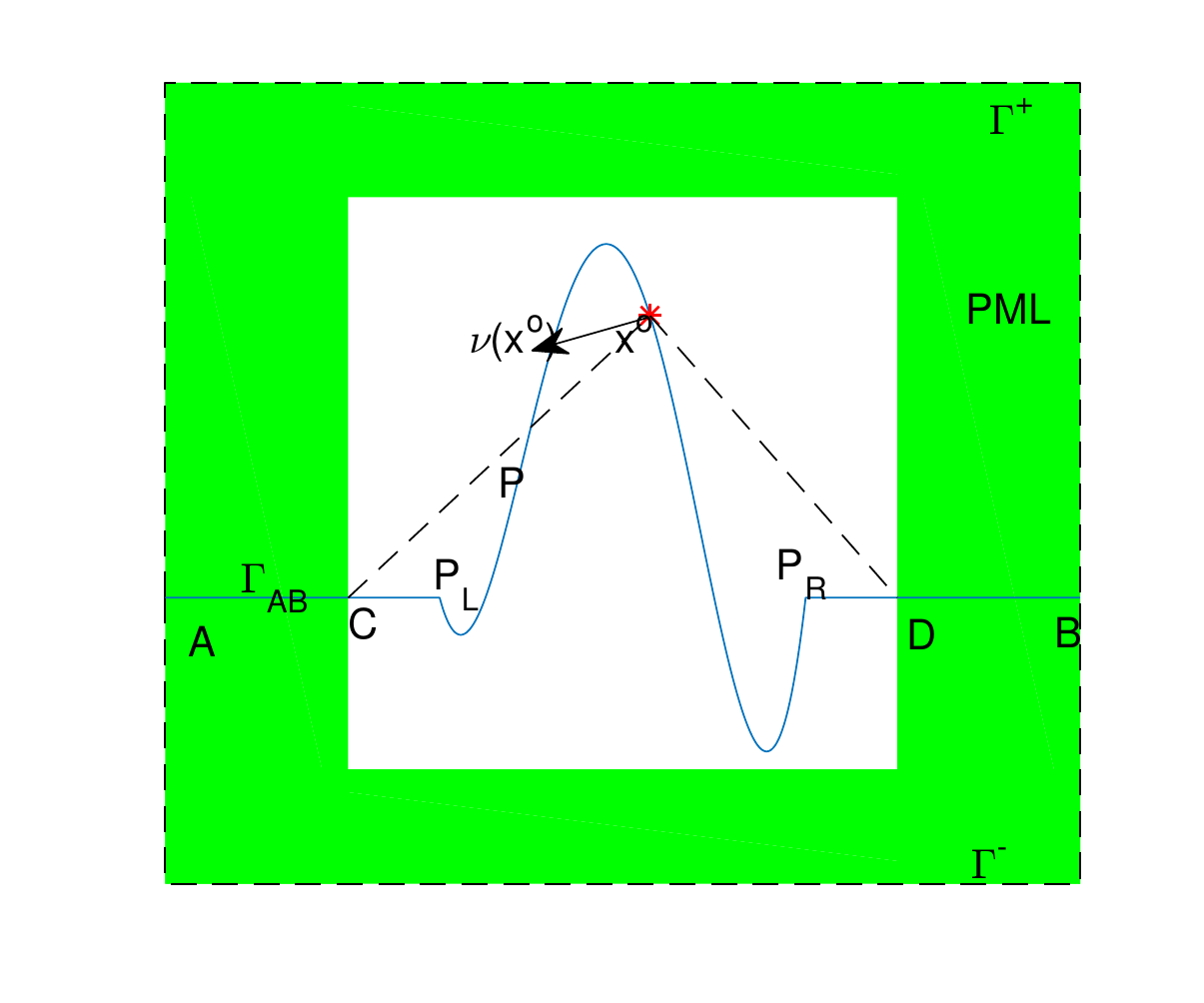}
  \caption{Profile of a 2D layered medium enclosed by a PML.}
		\label{fig:profile}
\end{figure}

A typical profile of a two-layer medium enclosed by a PML is shown in
Figure~\ref{fig:profile}. 
We now derive the Green's representation theorem for $\tilde{u}$ in the
bounded domain $\Omega^s$ enclosed by the same curve
$\Gamma^s=\Gamma^+\cup\Gamma_{\rm AB}$.

For $x\in\Omega^s$, 
\begin{align}
	\label{eq:grerep}
	\tilde{u}^s(x) &= \int_{\Omega^s} \delta_x(y) \tilde{u}^s(y) dy\nonumber\\
							 &= \int_{\Omega^s} [-\nabla_y\cdot({\bf A}(y)
\nabla_y\tilde{G}(x,y)) - k_0^2 n^2 J(y) \tilde{G}(x,y)]\tilde{u}^s({\bf
y})dy\nonumber\\
&= \int_{\Omega^s}\{-\tilde{u}^s(y)\nabla_y\cdot({\bf
			A}(y)\nabla_y\tilde{G}(x,y))+\tilde{G}(x,y)\nabla_y\cdot({\bf
		A}(y)\nabla_y\tilde{u}^s(y))\}dy\nonumber\\
		&= \int_{\Gamma_s}\{\partial_{{\bm \nu}_c}\tilde{G}(x,y)\tilde{u}^s(y) -
            \tilde{G}(x,y)\partial_{{\bm \nu}_c}\tilde{u}^s(y)\}ds(y),
%&=\int_{\Gamma_{\rm PML}^+\cup\Gamma_{\rm AB}}\{ \partial_{{\bm \nu}}
%\tilde{G}(x,y)\alpha_1(y)\tilde{u}(y)-\tilde{G}(x,y)\alpha_1(y)\partial_{{\bm \nu}}\tilde{u}(y)\}ds(y).
\end{align}
where the last equality follows from the second Green's identity
\cite{mcl00}, and we defined the conormal direction ${\bm \nu}_c = {\bf
A}^{T}{\bm \nu}$ and so the conormal derivative $\partial_{{\bm \nu}_c} =
{\bm \nu}_c\cdot \nabla= {\bm \nu}\cdot {\bf A}\nabla $.  

Similarly, in the limit case when $k_0\rightarrow 0^+$, one obtains the
representation formula 
\begin{align}
	\label{eq:rep:lap}
	\tilde{u}_0(x) =
	\int_{\tilde{\Gamma}}\{\partial_{{\bm \nu}_c}\tilde{G}_0(x,y)\tilde{u}_0(y)-
	\tilde{G}_0(x,y)\partial_{{\bm \nu}_c}\tilde{u}_0(y)\}ds(y),
\end{align}
for the complexified Laplacian equation 
\begin{equation}
	\label{eq:cmp:lap}
	\tilde{\Delta} u_0(\tilde{x}) = \nabla\cdot({\bf A} \nabla \tilde{u}_0(x))
	= 0,\quad x\in\Omega^s.
\end{equation}
Correspondingly, the related fundamental solution becomes, 
\begin{equation}
	\label{eq:green:complap}
	\tilde{G}_0(x,y) = -\frac{1}{2\pi}\log\rho(\tilde{x},\tilde{y}).
\end{equation}

\subsection{PML-NtD operator}
Based on the two representation formulae (\ref{eq:grerep}) and
(\ref{eq:rep:lap}), we are ready to develop the Neumann-to-Dirichlet (NtD)
operator on $\Gamma^s$.   

Since $\tilde{u}_0=1$ solves (\ref{eq:cmp:lap}), we get from
(\ref{eq:rep:lap}) that
\begin{equation}
	\label{eq:lap:1}
	1 = \int_{\Gamma^s}\partial_{{\bm \nu}_c}\tilde{G}_0(x,y) ds(y),
\end{equation}
when $x\in\Omega^s$.  
Now, (\ref{eq:grerep})$-\tilde{u}^s(x)\times$(\ref{eq:lap:1}) yields
\begin{equation}
	\label{eq:der:ntd:1}
	0 = \int_{\Gamma^s} \{(\partial_{{\bm \nu}_c} \tilde{G}(x,y) \tilde{u}^s(y) -
		\partial_{{\bm \nu}_c}\tilde{G}_0(x,y) \tilde{u}^s(x)) -
\tilde{G}(x,y)\partial_{{\bm \nu}_c} \tilde{u}^s(y)\}ds(y),
\end{equation}
which has a weakly-singular kernel \cite{lulu14}. 

Therefore, when $x$ approaches an observation point
$x^o\in\Gamma^s$, one easily reproduces for
the bounded domain $\Omega^s$, $\tilde{u}^s$ and $\partial_{{\bm \nu}_c}
\tilde{u}^s$ satisfies on $\Gamma^s$
\begin{align}
  \label{eq:ntd:bie}
	\tilde{\cal K}[\tilde{u}^s](x^o) - \tilde{\cal K}_0 [1](x^o)
	\tilde{u}^s(x^o) = \tilde{\cal S} [\partial_{{\bm \nu}_c} \tilde{u}^s](x^o),
\end{align}
where we have defined the following boundary integral operators in a PML
environment,
\begin{align}
  \label{eq:def:S}
	\tilde{\cal S}[\phi](x^o) &= 2\int_{\Gamma^s} \tilde{G}(x^o,y) \phi(y) ds(y),\\
  \label{eq:def:K}
	\tilde{\cal K}[\phi](x^o) &= 2\cint_{\Gamma^s} \partial_{{\bm \nu}_c}
	\tilde{G}(x^o,y)
  \phi(y) ds(y),\\
  \label{eq:def:K0}
	\tilde{\cal K}_0[\phi](x^o) &= 2\cint_{\Gamma^s} \partial_{{\bm
	\nu}_c} \tilde{G}_0(x^o,y)
  \phi(y) ds(y).
\end{align}

Consequently, we get the PML-NtD operator $\tilde{\cal N}= (\tilde{\cal
K}-\tilde{\cal K}_0 1)^{-1} \tilde{\cal S}$ on $\Gamma^s$, which maps
$\partial_{{\bm \nu}_c} \tilde{u}^s$ to $\tilde{u}^s$ on $\Gamma^s$.  

\subsection{Truncating PML-NtD operator onto $\Gamma_{\rm AB}$}
Unlike the slowly decaying wave $u^s$, $\tilde{u}^s$ and $\partial_{{\bm
\nu}_c}\tilde{u}^s$ decay exponentially at infinity so that it is reasonable to
impose $\tilde{u}^s\approx 0$ and $\partial_{{\bm \nu}_c} \tilde{u}^s\approx
0$ on $\Gamma^+$; see \cite{ber94, ber96} and section 4.1 below.  Therefore,
operators $\tilde{\cal K}$ and $\tilde{\cal S}$ in (\ref{eq:ntd:bie}) can be
truncated and defined onto curve $\Gamma_{AB}$ only, that is, for
$x^o\in\Gamma_{AB}$, 
\begin{align}
	\label{eq:ntd:bie:tr}
	\tilde{\cal K}_{AB}[\tilde{u}^s](x^o) - \tilde{\cal K}_0 [1](x^o)
	\tilde{u}^s(x^o) &\approx \tilde{\cal S}_{AB} [\partial_{{\bm \nu}_c}
\tilde{u}^s](x^o),
	%\label{eq:dtn:bie:tr}
	%\tilde{\cal K}_{AB}'[\partial_{{\bm \nu}_c}\tilde{u}](x^o)+\tilde{\cal
%K}_0 [1](x^o)
	%\partial_{{\bm \nu}_c^*}\tilde{u}(x^o)  &\approx \tilde{\cal
%T}_{AB}[\tilde{u}](x^o) +
	%\tilde{\cal S}_0'[1](x^o) \partial_{\tau^*}\tilde{u}(x^o),
\end{align}
where the definition of $\tilde{\cal K}_{AB}$ is the same as $\tilde{\cal
K}$ in (\ref{eq:def:K}) but with the integral domain replaced with
$\Gamma_{AB}$, etc.    

However, the integral $\tilde{\cal K}_0[1]$ cannot be truncated onto
$\Gamma_{AB}$ since the density function is nonzero on $\Gamma^+$.
Nevertheless, it turns out that 
\begin{align}
	\label{eq:K_0:1}
	\tilde{\cal K}_0[1](x^o) &= -\frac{\theta^*}{\pi},
	%\label{eq:S_0:1}
	%\tilde{\cal S}_0'[1](x^o) &= 0,
\end{align}
where $\theta^*$ is the interior angle of $x^o$ on $\Gamma_{\rm AB}$
{even when $x^o$ is in the PML; the proof will be
shown in the Appendix.  }  { Unfortunately, such a formula
	cannot be directly used near corners of $\Gamma_{\rm AB}$ since numerical
	discrepancies would appear there \cite{lulu14}.  
}   We now discuss how to remove the integral domain $\Gamma^+$ for operator $\tilde{\cal K}_0$.   

We distinguish two cases:  
\begin{itemize}
  \item[(1).]{Suppose $x^o\in \Gamma_{CD}$.  
As shown
in Figure~\ref{fig:profile}, for the closed curve 
\[
	\Gamma_{s}=\Gamma^+\cup AC\cup Cx^o\cup x^oD\cup DB,
\]
using (\ref{eq:K_0:1}), we see that
\begin{equation}
  \label{eq:Gs:K0}
  2\cint_{{\Gamma}_s} \partial_{{\bm \nu}_c} \tilde{G}_0(x^o,y) 1 ds(y) = -\angle
  Cx^oD/\pi,  
\end{equation}
where we note that $\angle Cx^oD$ denotes the interior angle. On the other
hand, one easily sees that 
\begin{equation}
  2\cint_{Cx^o\cup x^oD} \partial_{{\bm \nu}_c} \tilde{G}_0(x^o,y) 1 ds(y) = 0,  
\end{equation}
so that
\begin{equation}
  2\int_{\Gamma^+\cup AC\cup DB} \partial_{{\bm \nu}_c}
  \tilde{G}_0(x^o,y) 1 ds(y) = -\angle Cx^oD/\pi,  
\end{equation}
where the integral in fact becomes a Riemann integral. This implies that
\begin{align}
  \label{eq:trunc:K0:1}
	\tilde{\cal K}_0[1](x^o) &= -\angle Cx^oD/\pi + \tilde{\cal
  K}_{0,CD}[1](x^o),\nonumber\\
	&=-\angle Cx^oD/\pi + {\cal K}_{0,CD}[1](x^o),
\end{align}
where the subscript $CD$ indicates that the integral domain is
$\Gamma_{CD}$, and the second equality holds since the integral domain is
outside the PML.  Furthermore, since the integrated domain is outside the
PML, one easily gets \cite{lulu14}
\begin{equation}
  \label{eq:trunc:K0:12}
  \tilde{\cal K}_0[1](x^o) = -\angle Ax^oB/\pi + {\cal K}_{0,AB}[1](x^o).
\end{equation}
We remark that $(\ref{eq:trunc:K0:12})$ is more stable than
(\ref{eq:trunc:K0:1}) since $x^o$ is sufficiently far away from $A$ and $B$.
}
\item[(2).]{Suppose $x^o\in AC\cup DB$. Since now $x^o$ corresponds to a
	smooth point of $\Gamma_{AB}$ and since it is sufficiently far away from
	potential corners of $\Gamma_{CD}$, we can directly use the exact formula
	(\ref{eq:K_0:1})
	that $\tilde{\cal K}_0[1](x^o)=-1.$
}
\end{itemize}

After the truncation, the BIE (\ref{eq:ntd:bie:tr}) only depends on the bounded
curve $\Gamma_{AB}$. Therefore, by properly discretizing $\tilde{\cal
K}_{AB}$, $\tilde{\cal S}_{AB}$, and ${\cal K}_{0,AB}$, we are able to
approximate the PML-NtD operator $\tilde{\cal N}$ on $\Gamma_{AB}$ now.  

\section{Numerical implementation}
Suppose the piecewise smooth and open curve $\Gamma_{AB}$ is parameterized by
$x(s)=\{(x_1(s), x_2(s))|0\leq s\leq L\}$, where $s$ is the arclength.
Since $\Gamma_{AB}$ possibly contains corners, to smoothen the
non-differentiable $x(s)$, we construct a scaling function $s=w(t), 0\leq
t\leq 1$ following \cite{colkre13}, whose derivatives vanish at corners up
to order $p$. For example, for a smooth segement of $\Gamma_{AB}$
correpsonding to $s\in[s^0, s^1]$ and $t\in[t^0, t^1]$ such that
$s^l=w(t^l)$ for $l=0,1$ correspond to two corners, we may
take
\begin{equation}
	\label{eq:gfun}
	s=w(t) = \frac{s^0w_1^p + s^1w_2^p}{w_1^p+w_2^p},\quad t\in[t^0,t^1],
\end{equation}
where
\[
	w_1=\left(\frac{1}{2}-\frac{1}{p}\xi^3\right)+\frac{\xi}{p}+\frac{1}{2},\quad
	w_2 = 1-w_1,\quad \xi = \frac{2 t - (t^0+t^1)}{t^1-t^0}.
\]
Assume that $t\in[0,1]$ is uniformly sampled by an even number, denoted by
$N$, of grid points $\{t_j=jh\}_{j=1}^{N}$ with grid size $h=1/N$, and that
the grid points contain those corner points. The scaling function $s=w(t)$
creates a graded mesh on $\Gamma_{AB}$ in the sense that it makes part of
grid points cluster around corners while keeping the other part almost
uniformly spaced \cite{colkre13}.  

We shall discuss numerically discretizing the integral operators $\tilde{\cal
K}_{AB}$, $\tilde{\cal S}_{AB}$, and ${\cal K}_{0,AB}$ on $\Gamma_{AB}$ in
this section. To simplify the notations, we use $x(t)$ to denote
$x(w(t))$, and use $x'(t)$ to denote $\frac{dx}{ds}(w(t))w'(t)$.

\subsection{Setup of the PML}
Once reparameterized by parameter $t$, $x(t)$ now becomes at least a
$C^p$-class function and we can expect that integrands in
(\ref{eq:ntd:bie:tr}) are sufficiently smoothened near corners.
However, as $\Gamma_{AB}$ overlaps with the PML, if $\sigma_1$ in
(\ref{eq:sigma}) is not properly chosen, those integrands could have weaker
regularities at the entrance points $C$ and $D$, as shown in
Figure~\ref{fig:profile}, since $\tilde{x}_1(x_1)$ may not be smooth
there. We remark here that we do not need to specify $\sigma_2$ since
$\Gamma_{AB}$ is far away from the horizontal PML regions parallel to $x_1$.   

To ensure that $\tilde{x}_1(x_1)$ is at least a $C^p$-class function like
$x_1(t)$ at $C$ and $D$, we require that derivatives of $\sigma_1$ vanishes
at the entrances up to order $\tilde{p}$; to be on the safe side, we choose
$\tilde{p}=p+2$.   This motivates us to use a function similar to the
scaling function $w$ in (\ref{eq:gfun}) to construct $\sigma_1$. Suppose
the PML on $\Gamma$ is defined by $\{(x_1,0)| a_1\leq |x_1|\leq a_1+T\}$
where $T$ denotes the thickness of the PML.  Then, for $x_1\in[a_1,a_1+T]$,
we take
\begin{equation}
	\label{eq:sigma1}
	\sigma_1(x_1) = \frac{0 \tilde{w}_1^{\tilde{p}} +
		2ST\tilde{w}_2^{\tilde{p}}}{\tilde{w}_1^{\tilde{p}}+\tilde{w}_2^{\tilde{p}}},\quad
        x_1\in[a_1,a_1+T],
\end{equation}
where
\[
	\tilde{w}_1=\left(\frac{1}{2}-\frac{1}{\tilde{p}}\xi^3\right)+\frac{\xi}{\tilde{p}}+\frac{1}{2},\quad
	\tilde{w}_2 = 1-\tilde{w}_1,\quad \xi = \frac{2 x_1 -
        (a_1+T)}{T}.
\]
It is not hard to show that $\sigma_1$ bijectively maps $[a_1,a_1+T]$ to
$[0,ST]$, and satisfies the desired property at $x=a_1$. When
$x_1\in[-a_1-T,-a_1]$, one defines $\sigma_1(x_1)=\sigma_1(-x_1)$.    

Now we show how the PML absorbs an outgoing wave. Consider on
$\{(x_1,0)|x_1>a_1\}$, a simple outgoing wave $f^{s}(x)=e^{i c x_1}$ for a
given $c>0$ as $x_1\rightarrow\infty$. In the PML, we obtain
\begin{equation}
	\label{eq:us:pml}
	\tilde{f}^s(x)=f^s(\tilde{x}) = e^{ic\tilde{x}_1} = e^{icx_1}
	e^{-c\ {\rm imag}(\tilde{x}_1)},
\end{equation}
where 
\[
	{\rm imag}(\tilde{x}_1) = \int_{a_1}^{x_1}\sigma_1(t)dt.  
\]
Clearly, a larger $S$ produces a larger $\sigma_1$ so that the imaginary part
of $\tilde{x}_1$ becomes larger, and therefore $\tilde{f}^s(x)$ decays more
quickly and is absorbed more completely at the boundary $x_1=a_1+T$ of the PML.
We will refer to $S$ as the absorbing magnitude of the PML in the following.  

On the other hand, effectiveness of the PML is also closely related to
the magnitude of $c$; the greater $c$ is, the more effective the PML
becomes.    In general, our unknown outgoing wave $u^s(x)$ restricted on
$\Gamma$ contains many such simple outgoing functions $f^s(x)$ but with
different values of $c$.  One way to increase the smallest value of $c$
among those simple outgoing functions, is to place the PML sufficiently far
away from the local perturbation curve $P$. Empirically, for a plane
incident wave, the distance between the PML and curve $P$ can be around one
wavelength; for a cylindrical wave due to a source $x^*$, it is safer to place the
PML at lease one wavelength away from curve $P$ as well as the point source
$x^*$.

\subsection{Discretizing $\tilde{\cal S}_{AB}$}
According to its definition, $\tilde{\cal S}_{AB}$ acting on
$\partial_{{\bm \nu}_c}\tilde{u}^s$ at $x=x(t_l)$, $l=1,\cdots, N$ can be
parameterized by
\begin{equation}
	\label{eq:tS:AB}
	\tilde{\cal S}_{AB}[\partial_{{\bm \nu}_c}\tilde{u}^s](x(t_l)) = \int_0^1
	S(t_l,t) \phi(t) dt,  
\end{equation}
where
\begin{align}
	\label{eq:def:mS}
	S(t_l,t) &= \frac{i}{2} H_0^{(1)}(k_0n\ \d(t_l,t)),\\
	\d(t_l,t) &= \rho(x(t_l),x(t)),\\
	\phi(t) &= \partial_{{\bm \nu}_c}\tilde{u}^s(x(t))|x'(t)|.
\end{align}
Clearly, $\partial_{{\bm \nu}_c}\tilde{u}^s$ is not continuous on $\Gamma_{AB}$
since ${\bm \nu}_c$	is discontinuous at corners. However, since $x'(t)$
vanishes at corners, the scaled co-normal derivative $\phi(t)$ is
smoothened.  

One way to discretize the integral in (\ref{eq:tS:AB}) is using the kernel
splitting technique developed in \cite{colkre13}. Specifically, the
logarithmic singularity of $S$ at $t=t_l$ can be splitted out in terms of
\[
	S(t_l,t)=S_1(t_l,t)\log(4\sin^2(\pi(t_l-t)) + S_2(t_l,t),
\]
where for $t\neq t_l$, we have
\begin{equation}
	\label{eq:S1}
	S_1(t_l,t) = -\frac{1}{2\pi}J_0(k_0 n\ \rho(x(t_l), x(t))).
\end{equation}

However, such a technique loses high accuracy when the argument of $J_0$ in
(\ref{eq:S1}) becomes complex. Specifically, if $x(t_l)$ or $x(t)$ is in the
PML, $\rho(x(t_l), x(t))$ may have large imaginary part, giving rise to a
blow up function $J_0$ and hence inducing numerical instabilities.  

Fortunately, to treat integrands with logarithmic singularities, Alpert
\cite{alp99} developed an efficient quadature rule which does not require a
kernel splitting process. Following such an approach, we may discretize the
integral in (\ref{eq:tS:AB}) as
\begin{align}
	\label{eq:dis:tS}
	\tilde{\cal S}_{AB}[\partial_{{\bm \nu}_c}\tilde{u}^s](x(t_l)) \approx
	\sum_{k=1}^{K_1}& \gamma_k h[S(t_l, t_l+\delta_kh) \phi(t_l+\delta_k h)
		\nonumber\\
		&+ S(t_l, t_l + 1-\delta_kh) \phi(t_l+1-\delta_k h)] \nonumber\\
		+ \sum_{k=K_2}^{N-K_2}& hS(t_l, t_l+t_k) \phi(t_l+t_k)\nonumber\\
	=\sum_{k=1}^{K_1}& \gamma_k h[S(t_l, t_l+\delta_kh) \phi(t_l+\delta_k h)
		\nonumber\\
		&+ S(t_l, t_l -\delta_kh) \phi(t_l-\delta_k h)] \nonumber\\
	+ \sum_{k=K_2}^{N-K_2}& hS(t_l, t_{{\rm mod}(l+k,N)}) \phi(t_{{\rm
            mod}(l+k,N)}),
\end{align}
where values of $K_1$, $K_2$, $\gamma_k$, and $\delta_k$ depend on the order
of Alpert's quadrature rule and can be precomputed. For example, in a $6$-th
order quadrature formula, we have $K_1=5$ and $K_2=3$; the associated
$\{\delta_k,\gamma_k\}_{k=1}^{5}$ are given in Table~\ref{table:alpt:6th}.
\begin{table}[h]
  \centering
  \begin{tabular}{c|c|c}
    $k$ & $\delta_k$ &  $\gamma_k$ \\\hline
    1 & 4.00488 41949 26570 E-03 & 1.67187 96911 47102 E-02 \\\hline
    2 & 7.74565 53733 36686 E-02 & 1.63695 83714 47360 E-01 \\\hline
    3 & 3.97284 99935 23248 E-01 & 4.98185 65697 70637 E-01 \\\hline
    4 & 1.07567 33529 15104 E$+$00 & 8.37226 62455 78912 E-01 \\\hline
    5 & 2.00379 69271 11872 E$+$00 & 9.84173 08440 88381 E-01 \\\hline
  \end{tabular}
  \caption{The $6$-th order Alpert's quadrature rule.  }
  \label{table:alpt:6th}
\end{table}

On the other hand, for sufficiently large $p$, it is reasonable to regard
$\phi(t)$ as a smooth periodic function so that we may approximate $\phi$ by
its trigonometric interpolation \cite{tre00}
\begin{align}
	\label{eq:f:trg}
	\phi(t)\approx\sum_{j=1}^{N} \phi(t_j)L(t-t_j),
\end{align}
where $L(t)=\sin(N\pi t)/[N\tan(\pi t)]$ is the Sinc function, satisfying
$L(t_j)=0$ for $1\leq j< N$ and $L(1)=L(0)=1$. Utilizing (\ref{eq:f:trg}),
we may rewrite equation (\ref{eq:dis:tS}) in terms of $\phi(t_j)$ for $1\leq
j\leq N$ so that we obtain an $N\times N$ matrix ${\bf S}$ that satisfies
\begin{equation}
	\tilde{\cal S}_{AB}[\partial_{\nu_c}\tilde{u}^s]\left[
		\begin{array}{c}
			x(t_1)\\
			\vdots\\
			x(t_N)
		\end{array}
	\right] \approx {\bf S} \left[
		\begin{array}{c}
			\phi(t_1)\\
			\vdots\\
			\phi(t_N)
		\end{array}
	\right],
\end{equation}
where the left-hand side represents a column vector produced by evaluating
$\tilde{\cal S}_{AB}[\partial_{\nu_c}\tilde{u}^s]$ at each
element $x(t_j)$ of the column vector for $1\leq j\leq N$.    

\subsection{Discretizing $\tilde{\cal K}_{AB}$}
According to its definition, $\tilde{\cal K}_{AB}$ acting on $\tilde{u}^s$ at 
$x=x(t_l)$ can be parameterized as
\begin{equation}
	\label{eq:tK:AB}
	\tilde{\cal K}_{AB}[\tilde{u}^s](x(t_l)) = \int_{0}^{1} K(t_l,t)g(t)dt,
\end{equation}
where 
\begin{align}
	\label{eq:def:tKAB}
	K(t_l,t) &= -\frac{ik_0n}{2}\frac{\kappa(t_l,t)}{\d(t_l,t)} H_1^{(1)}(k_0
	n\d(t_l,t)),\\
	\label{eq:def:kappa}
	\kappa(t_l,t) &= \tilde{x}_2'(t)(\tilde{x}_1(t) - \tilde{x}_1(t_l))
	-\tilde{x}_1'(t)(\tilde{x}_2(t)-\tilde{x}_2(t_l)),\\
	g(t) &= \tilde{u}^s(x(t)).
\end{align}
Thus, similar to operator $\tilde{\cal S}_{AB}$, by appling Alpert's quadrature
rule, we discretize the integral in (\ref{eq:tK:AB}) as
\begin{align}
	\label{eq:dis:tK}
	\tilde{\cal K}_{AB}[\tilde{u}](x(t_l)) \approx \sum_{k=1}^{K_1}& \gamma_k
	h[K(t_l, t_l+\delta_kh) g(t_l+\delta_k h) \nonumber\\ &+ K(t_l,
t_l -\delta_kh) g(t_l-\delta_k h)] \nonumber\\ 
+ \sum_{k=K_2}^{N-K_2}& hK(t_l, t_{{\rm mod}(l+k,N)}) g(t_{{\rm
mod}(l+k,N)}).
\end{align}

By choosing $p$ sufficiently large, we may approximate $g$ by its
trigonometric interpolation like (\ref{eq:f:trg}) but with $\phi$ replaced by
$g$. Consequently, we may rewrite equation (\ref{eq:dis:tK}) in terms of
$g(t_j)$ for $1\leq j\leq N$ 
so that we obtain an $N\times N$ matrix ${\bf K}$ that satisfies
\begin{equation}
	\tilde{\cal K}_{AB}[\tilde{u}^s]\left[
		\begin{array}{c}
			x(t_1)\\
			\vdots\\
			x(t_N)
		\end{array}
	\right] \approx {\bf K} \left[
		\begin{array}{c}
			g(t_1)\\
			\vdots\\
			g(t_N)
		\end{array}
	\right].
\end{equation}
The discretization of ${\cal K}_{0, AB}$ in (\ref{eq:trunc:K0:12}) can be
derived similarly. 

After the discretization of $\tilde{\cal S}_{AB}$, $\tilde{\cal K}_{AB}$,
and ${\cal K}_{0, AB}$, one obtains from (\ref{eq:ntd:bie:tr}) that
\begin{equation}
	\label{eq:bie:mat}
	({\bf K} - {\bf H}) \tilde{\bf u}^s  \approx {\bf S} {\bm \phi},
\end{equation}
where ${\bf H}$ is a diagonal matrix with entries $\tilde{\cal
K}_0[1](x(t_l))$ for $1\leq l\leq N$, 
\begin{align*}
\tilde{\bf
u}^s&=[g(x(t_1)),\ldots, g(x(t_N))]^{T},\\
{\bm
\phi}&=[{\phi}(x_(t_l)),\ldots, {\phi}(x(t_N))]^{T}.
\end{align*}

Consequently, one gets
\begin{equation}
	\label{eq:ntd:mat}
	\tilde{\bf u}^s \approx ({\bf K}-{\bf H})^{-1}{\bf S} {\bm \phi}:= {\bf N}
	{\bm \phi},
\end{equation}
where the $N\times N$ matrix ${\bf N}$ in fact approximates a scaled PML-NtD
operator $\tilde{\cal N}_s$ which maps $\phi = \partial_{{\bm \nu}_c}
\tilde{u}^s(x)|x'|$ to $\tilde{u}^s$ on $\Gamma_{AB}$.  

\subsection{A stabilizing technique}
Clearly, to make the approximations of $\tilde{\cal S}_{AB}$ and
$\tilde{\cal K}_{AB}$ accurate enough, a high order quadrature rule and a
large scaling parameter $p$ are always preferable; otherwise, one needs a
sufficiently large $N$. Suppose we desire a $6$-th order of accuracy so that
nodes and weights of Alpert's quadrature rule are chosen based on
Table~\ref{table:alpt:6th}. To be consistent, we choose
$p=6$ in the scaling function $s=w(t)$.  Under such a circumstance, when
computing the kernel functions $S(t_l, t)$ and $K(t_l, t)$, we observe that
$|t_l-t|$ can be as small as $\delta_1 h = O(\frac{10^{-3}}{N})$.  When
$t_l$ is close to a corner point, the physical distance $\d(t_l, t)$ can be
further shrunk to $O(\frac{10^{-3p}}{N^p})=O(\frac{10^{-18}}{N^6})$ by the
scaling function.  Unfortunately, even for a coarse mesh, this can
be less than or close to the round-off error $O(10^{-16}x(t_l))$ in the
computation of $\d(t_l,t)$. In such a situation, $\d(t_l,t)$
is simply regarded as $0$ in a double-precision computation.
Consequently, division by zero occurs in the computation of $S(t_l, t)$ and
$K(t_l,t)$ when $t$ is close to $t_l$ and when $t_l$ is close to a corner.  

To resolve this instability issue, one approach is to reduce the accuracy
order to $p=3$ or less. However, this can make the total computational
process extremely inefficient in practice. In this section, we develop
numerical techiques which can accurately compute $S(t_l,t)$ and $K(t_l, t)$
when $t$ is close to $t_l$ and $t_l$ is close to a corner.  

Observing their expressions (\ref{eq:def:mS}) and (\ref{eq:def:kappa}),
the instability issue comes from the two terms $\d(t_l, t)$ and $\kappa(t_l, t)$
since they involve subtractions of two extremely close quantities.  
We discuss $\d(t_l, t)$ first. 

Without loss of generality, we assume $t>t_l$, so that
$\tilde{x}(\xi), \xi\in[t_l,t]$ becomes a piecewise smooth function; note
that here $\tilde{x}(\xi)$ may meet some corner. At first, we assume that
$\tilde{x}$ on $[t_l, t]$ is smooth. 

To preserve enough significant digits, we require accurately computing 
\begin{equation}
  \label{eq:diffxi}
  \tilde{x}_i(t)-\tilde{x}_i(t_l), 
\end{equation}
for $i=1,2$.  There are two approaches to realize this. The first approach
is to use the Taylor series of $\tilde{x}_i$ at $t_l$, that is,   
\begin{align}
  \label{eq:taylor}
  \tilde{x}_i(t)-\tilde{x}_i(t_l) = \sum_{j=1}^{\infty}
  \frac{\tilde{x}_i^{(j)}(t_l)}{j!}(t-t_l)^{j}.
\end{align}
Unfortunately, this approach is ineffective since it is not easy to control
the truncation error and since we require the computation of many high order
derivatives. The second and more effective approach utilizes the
Newton-Lebnitz formula, rewriting (\ref{eq:diffxi}) in the form,
%\begin{align}
%  \label{eq:newton1}
%  \tilde{x}_i(t)-\tilde{x}_i(t_l) = \int_{t_l}^{t}
%  \frac{d\tilde{x}_i}{ds}(w(\tau))w'(\tau) d\tau,
%\end{align}
%and 
\begin{align}
  \label{eq:newton2}
  \tilde{x}_i(t)-\tilde{x}_i(t_l) &= \tilde{x}_i(w(t)) - \tilde{x}_i(w(t_l)) \nonumber\\
  &= \int_{0}^{w(t)-w(t_l)} \frac{d \tilde{x}_i}{ds}(w(t_l)+s)ds,\nonumber\\
  &=\int_{0}^{\int_{t_l}^{t}w'(\tau)d\tau}\frac{d \tilde{x}_i}{ds}(w(t_l)+s)ds.  
\end{align}
for $i=1,2$. Such an representation gives rise to siginficant advantages.
Specifically, the integrand in the primary integral is an $O(1)$
quantity so that using numerical integrations (e.g., Gaussian quadrature
rules), to compute the integral can highly reduce round-off errors;
moreover, we only require the first-order derivative of $\tilde{x}_i$ to
obtain accurate result.  To ensure stability, the upper limit is also
rewritten as an integral form.  Consequently, $\d(t_l,t)$ can be evaluated via
\begin{align}
  \label{eq:eval:dist}
  \d(t_l, t) =
  \sqrt{\sum_{i=1}^2\left(\int_{0}^{\int_{t_l}^{t}w'(\tau)d\tau}\frac{d
  \tilde{x}_i}{ds}(w(t_l)+s)ds\right)^2}.
\end{align}
We remark that the aim of using arclength $s$ but not the grading parameter
$t$ as the integral variable is to further stabilize the involved
computations since integrands roughly become $O(1)$ quantities.

Next, we discuss the computation of 
\begin{align}
  \label{eq:kappa:2}
  \kappa(t_l,t)=&
  w'(t)\Big[\frac{d\tilde{x}_2}{ds}(w(t))\left(\tilde{x}_1(w(t))-\tilde{x}_1(w(t_l))\right)\nonumber\\
  &-\frac{d\tilde{x}_1}{ds}(w(t))\left(\tilde{x}_2(w(t))-\tilde{x}_2(w(t_l))\right)\Big]\nonumber\\
:=& w'(t)\bar{\kappa}(t_l,t).
\end{align}
Using Newton-Lebnitz formula, we may rewrite $\bar{\kappa}(t_l,t)$ as 
\begin{align}
  \label{eq:eval:kappa}
  \bar{\kappa}(t_l,t) =&\int_{0}^{w(t)-w(t_l)}\Big[\frac{d^2\tilde{x}_2}{ds^2}(w(t_l)+s)\left(\tilde{x}_1(w(t_l)+s)-\tilde{x}_1(w(t_l))\right)\nonumber\\
  &-\frac{d^2\tilde{x}_1}{ds^2}(w(t_l)+s)\left(\tilde{x}_2(w(t_l)+s)-\tilde{x}_2(w(t_l))\right)\Big]ds\nonumber\\
  =&\int_{0}^{\int_{t_l}^{t}w'(\tau)d\tau}\int_{0}^{s}\Big[\frac{d^2\tilde{x}_2}{ds^2}(w(t_l)+s)\frac{d\tilde{x}_1}{ds}(w(t_l)+\eta)\nonumber\\
  &-\frac{d^2\tilde{x}_1}{ds^2}(w(t_l)+s)\frac{d\tilde{x}_2}{ds}(w(t_l)+\eta)\Big]d\eta ds.  
\end{align}
Using numerical integrations to compute the above double integrals can yield
accurate results.   

Now, suppose that $\tilde{x}(\xi), \xi\in[t_l,t]$ contains a corner at
$\xi=t^*\in(t_l,t)$. Since $\tilde{x}(\xi)$ consists of two smooth segments
corresponding to $[t_l,t^*]$ and $[t^*,t]$, respectively, the following
splitting
\begin{equation}
  \label{eq:diffx_i:2}
  \tilde{x}_i(t)-\tilde{x}_i(t_l) = (\tilde{x}_i(t)-\tilde{x}_i(t^*)) +
  (\tilde{x}_i(t^*)-\tilde{x}_i(t_l)),
\end{equation}
indicates that Newton-Lebnitz formula is applicable for either term on the
right-hand side so that numerical integrations can still offer an accurate
result for $\d(t_l, t)$. 

As for $\kappa(t_l, t)$, we have
\begin{align}
  \label{eq:barkappa}
  \bar{\kappa}(t_l, t) =& \bar{\kappa}(t_l,t^*) + \bar{\kappa}(t^*,t)
  \nonumber\\
  &+
  \Big[\left(\frac{d\tilde{x}_2}{ds}(w(t))-\frac{d\tilde{x}_2}{ds}(w(t^*+))\right)\left(\tilde{x}_1(w(t^*))-\tilde{x}_1(w(t_l))\right)\nonumber\\
  &-\left(\frac{d\tilde{x}_1}{ds}(w(t))-\frac{d\tilde{x}_1}{ds}(w(t^*+))\right)\left(\tilde{x}_2(w(t^*))-\tilde{x}_2(w(t_l))\right)\Big]\nonumber\\
  &+\Big[\left(\frac{d\tilde{x}_2}{ds}(w(t^*+))-\frac{d\tilde{x}_2}{ds}(w(t^*-))\right)\left(\tilde{x}_1(w(t^*))-\tilde{x}_1(w(t_l))\right)\nonumber\\
  &-\left(\frac{d\tilde{x}_1}{ds}(w(t^*+))-\frac{d\tilde{x}_1}{ds}(w(t^*-))\right)\left(\tilde{x}_2(w(t^*))-\tilde{x}_2(w(t_l))\right)\Big]\nonumber\\
  =& \bar{\kappa}(t_l,t^*) + \bar{\kappa}(t^*,t) \nonumber\\
  &+\int_{0}^{\int_{t^*}^{t}w'(\tau)d\tau}\int_{0}^{\int_{t_l}^{t^*}w'(\tau)d\tau}\Big[\frac{d^2\tilde{x}_2}{ds^2}(w(t^*)+s)\frac{d\tilde{x}_1}{ds}(w(t_l)+\eta)\nonumber\\
  &-\frac{d^2\tilde{x}_1}{ds^2}(w(t^*)+s)\frac{d\tilde{x}_2}{ds}(w(t_l)+\eta)\Big]d\eta
    ds\nonumber\\
    &+\Big[\left(\frac{d\tilde{x}_2}{ds}(w(t^*+))-\frac{d\tilde{x}_2}{ds}(w(t^*-))\right)\int_{0}^{\int_{t_l}^{t^*}w'(\tau)d\tau}\frac{d\tilde{x}_1}{ds}(w(t_l)+\eta)
		d\eta\nonumber\\
  &-\left(\frac{d\tilde{x}_1}{ds}(w(t^*+))-\frac{d\tilde{x}_1}{ds}(w(t^*-))\right)\int_{0}^{\int_{t_l}^{t^*}w'(\tau)d\tau}\frac{d\tilde{x}_2}{ds}(w(t_l)+\eta)
		d\eta\Big],
\end{align}
where $(t^*+)$ and $(t^*-)$ indicate limits are taken from right side and
left side, respectively.  Clearly, all the four terms on the right-hand side
can be accurately evaluated through numerical integrations.   

\section{Wave field evaluations}
Suppose now in each domain $\Omega_j$, we have obtained an $N\times N$ matrices
${\bf N}_j$ to approximate the scaled NtD operator $\tilde{\cal N}_{s,j}$,
mapping $|x'|\partial_{{\bm \nu}_c} \tilde{u}^{s}_j$ to $\tilde{u}^s_j$ on
$\Gamma_{AB}$, for $j=1,2$.  Then
\begin{equation}
	\label{eq:sNtD}
	{\bf N}_j {\bm \phi}_j= \tilde{\bf u}_j^s ,
\end{equation}
where
\begin{align*}
	\tilde{\bf u}_j^s &= [\tilde{u}_j^s(x(t_1)),\ldots,
\tilde{u}_j^s(x(t_N))]^{T},\\
	{\bm \phi}_j &= [|x'(t_1)|\partial_{{\bm \nu}_c} \tilde{u}_j^s(x(t_1)),\ldots,
|x'(t_N)|\partial_{{\bm \nu}_c} \tilde{u}_j^s(x(t_N))]^{T}.
\end{align*}

According to the transmission conditions (\ref{eq:cond:us1}) and
(\ref{eq:cond:us2}), the complexified outgoing wave $\tilde{u}_j^s$, at
the $N$ grid points on $\Gamma_{AB}$, satisfies
\begin{align}
	\label{eq:cond:cus1}
	\tilde{\bf u}_1^s - \tilde{\bf u}_2^s &= {\bf b}_1,\\
	\eta_1{\bm \phi}_1 - \eta_2{\bm \phi}_2 &= {\bf b}_2,
\end{align}
where we have defined
\begin{align*}
	{\bf b}_1 &= [-[\tilde{u}_0^{tot}](x(t_1)),\ldots, -[\tilde{u}_0^{tot}](x(t_N))]^{T},\\
	{\bf b}_2 &= [-|x'(t_1)|[\eta_j\partial_{{\bm \nu}_c} \tilde{u}_0^{tot}](x(t_1)),\ldots,
	-|x'(t_N)|[\eta_j\partial_{{\bm \nu}_c} \tilde{u}_0^{tot}](x(t_N))]^{T}.
\end{align*}
Thus, by (\ref{eq:sNtD}), we obtain
\begin{align}
	\label{eq:linsys}
	\left[
\begin{array}{cc}
	{\bf N}_s^1 & - {\bf N}_s^2 \\
	\eta_1{\bf I} & -\eta_2{\bf I} 
\end{array}
	\right]\left[
		\begin{array}{c}
			{\bm \phi}_1\\
			{\bm \phi}_2
		\end{array}
	\right]
	= \left[
		\begin{array}{c}
			{\bf b}_1\\
			{\bf b}_2
		\end{array}
	\right],
\end{align}
which can be solved by
\begin{equation}
	\label{eq:sol:phi12}
\left[
		\begin{array}{c}
			{\bm \phi}_1\\
			{\bm \phi}_2
		\end{array}
	\right] = \left[
\begin{array}{cc}
	{\bf N}_s^1 & - {\bf N}_s^2 \\
	\eta_1{\bf I} & -\eta_2{\bf I} 
\end{array}
	\right]^{-1}\left[
		\begin{array}{c}
			{\bf b}_1\\
			{\bf b}_2
		\end{array}
	\right],
\end{equation}
or equivalently,
\begin{align}
	\label{eq:phi1:1}
	{\bm \phi}_1 &= ({\bf N}^1_s - \frac{\eta_1}{\eta_2}{\bf
	N}^2_s)^{-1}\left(\eta_2^{-1}{\bf N}^2_s{\bf b}_2+{\bf b}_1\right),\\
	{\bm \phi_2} &= \frac{\eta_1}{\eta_2}{\bm \phi}_1 - \frac{{\bf
	b}_2}{\eta_2}.
\end{align}
Consequently, we obtain $\tilde{\bf u}_j^s={\bf N}_j {\bm\phi}_j$ on
$\Gamma_{AB}$.

As for any point $x\in\Omega_j$, we may directly use (\ref{eq:grerep}) to
compute $\tilde{u}_j^s(x)$, that is
\begin{equation}
	\label{eq:field:evaluation}
	\tilde{u}_j^s(x) \approx \int_{\Gamma_{AB}}
	\{\partial_{{\bm
	\nu}_c}\tilde{G}_j(x,y)\tilde{u}_j^s(y)-\tilde{G}_j(x,y)\partial_{{\bm
\nu}_c} \tilde{u}_j^s(y)\} ds(y),
\end{equation}
where we keep curve $\Gamma_{AB}$ only since both
$\tilde{u}_j^s(y)$ and $\partial_{{\bm \nu}_c}\tilde{u}_j^s(y)$ approximately are $0$ on
$\tilde{\Gamma}/\Gamma_{AB}$. After parameterized by the scaling function $s=w(t)$
in (\ref{eq:gfun}), the integrand in (\ref{eq:field:evaluation}) becomes
periodic and smooth so that by trapezoidal rule, we may approximate
\begin{align}
	\label{eq:tildeu:tr}
	\tilde{u}_j^s(x^o) \approx \frac{1}{N}\sum_{l=1}^N&
	\Big[\partial_{{\bm \nu}_c}\tilde{G}_j^s(x^o,x(t_l))|x'(t_l)|\tilde{u}_j^s(x(t_l))
	\nonumber\\
	&-\tilde{G}_j(x^o,x(t_l)) |x'(t_l)|\partial_{{\bm \nu}_c}
	\tilde{u}_j^s(x(t_l))\Big]. 
\end{align}
Therefore, in the domain outside the PML, we obtain $u_j^s=\tilde{u}_j^s$
so that the total wave field $u^{tot}=u^s+u_0^{tot}$.  

\section{Numerical examples}
In this section, we will carry out several numerical experiments to
illustrate the proposed methodology. In all examples, the physical region is
defined as $\{(x_1,x_2)||x_1|\leq a_1, a_1>0\}$, while the PML region is
defined as $\{(x_1,x_2)|a_1\leq|x_1|\leq a_1+T, a_1>0, T>0\}$ with
thickness $T$. Therefore, the truncated interface $\Gamma_{AB}$ is just
$\Gamma$ restricted on $x_1\in[-a_1-T,a_1+T]$, while the physical region 
on $\Gamma_{AB}$, denoted by $\Gamma_P$ below, is just $\Gamma$ restricted
on $x_1\in[-a_1,a_1]$. To achieve a high-order accuracy, we take $p=6$ in
the scaling function $s=w(t)$ associated with the $6$-th order Alpert's
quadrature rule, using nodes and weights defined in
Table~\ref{table:alpt:6th}. We will mainly consider TM-polarization
problems.      

\subsection{Example 1: Perfectly flat surface}
To validate our method, the first example is a perfectly flat surface
$\Gamma=\{(x_1,x_2)|x_2=0\}$, where $n_1=1$, $n_2=2$, and the freespace
wavelength $\lambda=1$ so that $k_0=2\pi$.  When $u^{inc}$ represents a
plane incident wave, $u^{tot}=u_{0}^{tot}$ in (\ref{eq:u0}) is the exact
solution, making $u^s=0$ in both $\Omega_1$ and $\Omega_2$. To avoid
such trivial solutions, we here test the case when $u^{inc}$ is a
cylindrical wave due to a point source $x^*=(0,0.1)$, so that $u^{tot}$
represents a layered Green's function at $x^*$.  

In the implementation, although $\Gamma$ is smooth, we still set $(0,0)$ as
an artificial corner since it is close to the source $x^*$. As shown in
\cite{perarabru14}, an explicit expression of the layered Green's function is
available so that we can obtain the exact solution $u^{tot}_{exa}$ for
reference. 

\begin{figure}[!ht]
  \centering
  (a)\includegraphics[width=0.4\textwidth]{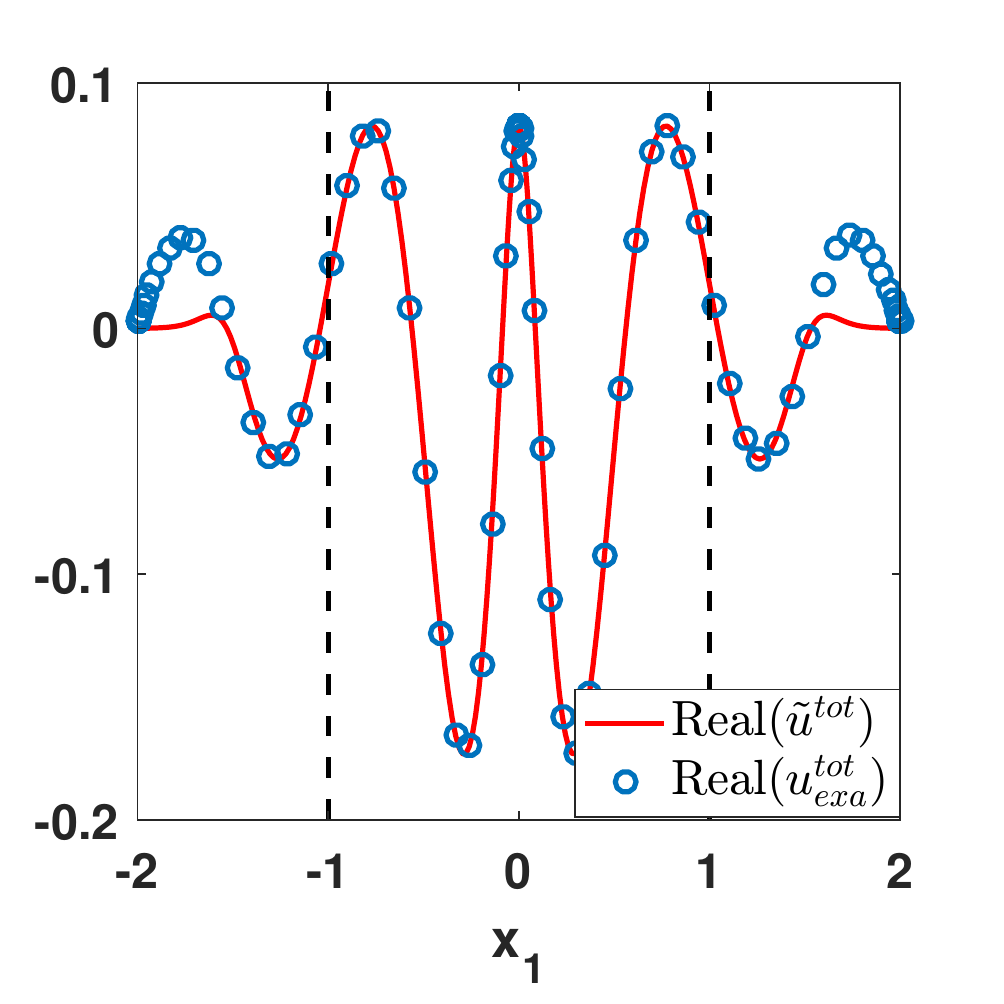}
  (b)\includegraphics[width=0.4\textwidth]{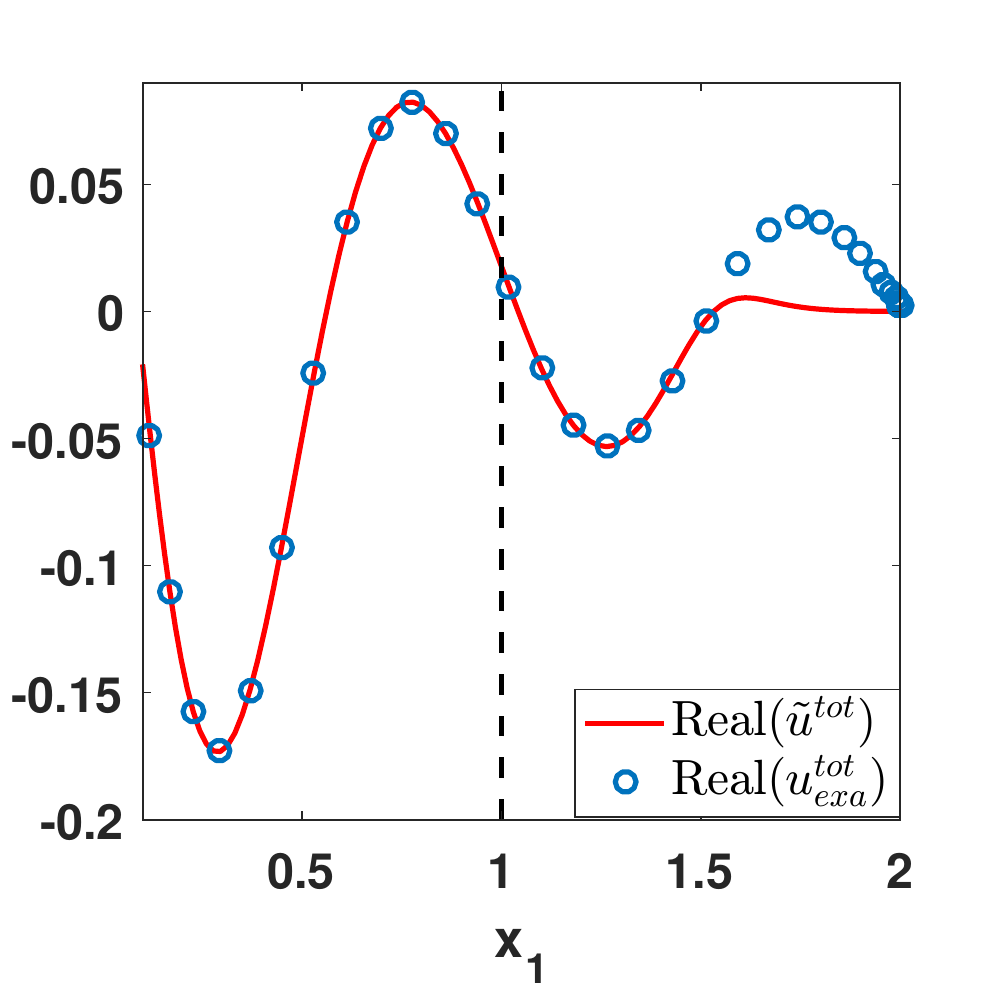}
  \caption{Example 1: real parts of $\tilde{u}^{tot}$ and the exact
    solution $u^{tot}_{exa}$ on: (a) $x_1\in[-1,1]$; (b) $x_1\in[0.3, 1]$.
	Dashed lines indicate entrances of the PML.  }
    \label{fig:ex1:cmp}
\end{figure}

\begin{figure}[!ht]
  \centering
  (a)\includegraphics[width=0.4\textwidth]{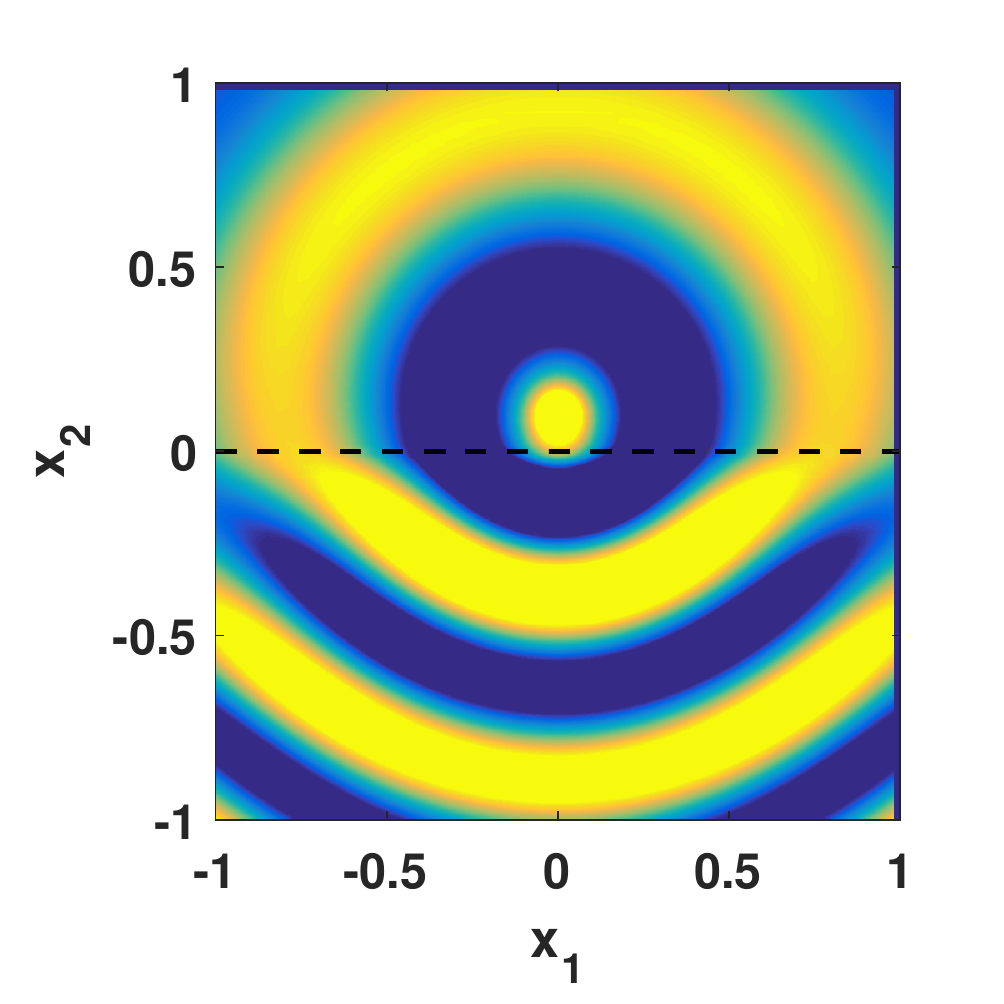}
  (b)\includegraphics[width=0.4\textwidth]{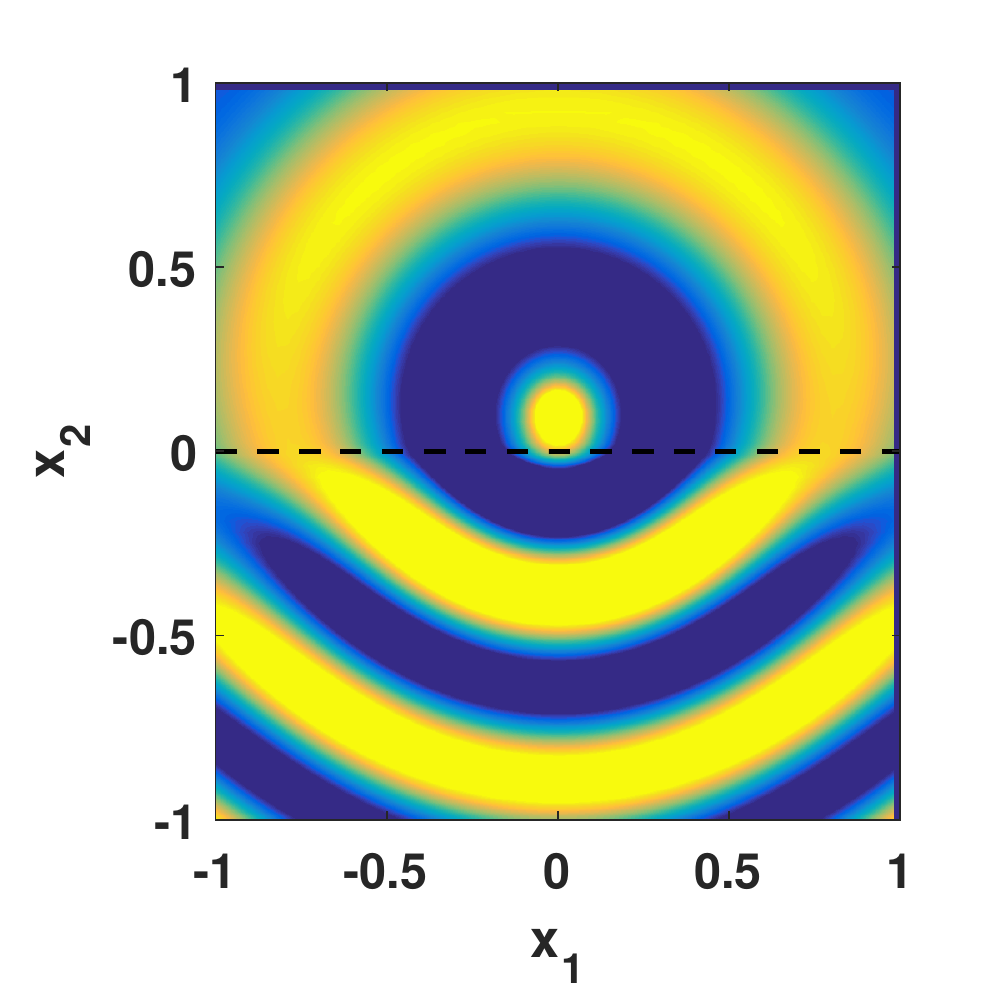}
	\caption{Example 1: real part of $u^{tot}$ on $[-1,1]\times[-1,1]$. (a)
		numerical solution; (b) exact solution.  Dashed line indicates location
	of $\Gamma$.  }
    \label{fig:ex1:cmp2}
\end{figure}

Taking $N=400$, $a_1=1$ and $T=1$, we compute $\tilde{u}^{tot}$, and
compare it with the exact solution $u^{tot}_{exa}$ on $\Gamma_{AB}$, as
shown in Figure~\ref{fig:ex1:cmp}.  
Clearly, on $\Gamma_P$, $\tilde{u}^{tot}$ and
$u^{tot}_{exa}$ coincide very well; in the PML region corresponding to
$|x_1|\in[1,2]$, $\tilde{u}^{tot}$ decay quickly to $0$ and $u^{tot}_{exa}$
still oscillates with a slowly decaying amplitude, as what we are expecting.
Figure~\ref{fig:ex1:cmp2}
show numerical and exact solutions of the real part of $u^{tot}$ in a box
$[-1,1]\times[-1,1]$, where Figure~\ref{fig:ex1:cmp2}(a) is based on a 
numerical solution using $N=400$ grid points on $\Gamma_{AB}$.  Obviously, they
coincide with each other quite well.   

\begin{figure}[!ht]
  \centering
	(a)\includegraphics[width=0.4\textwidth]{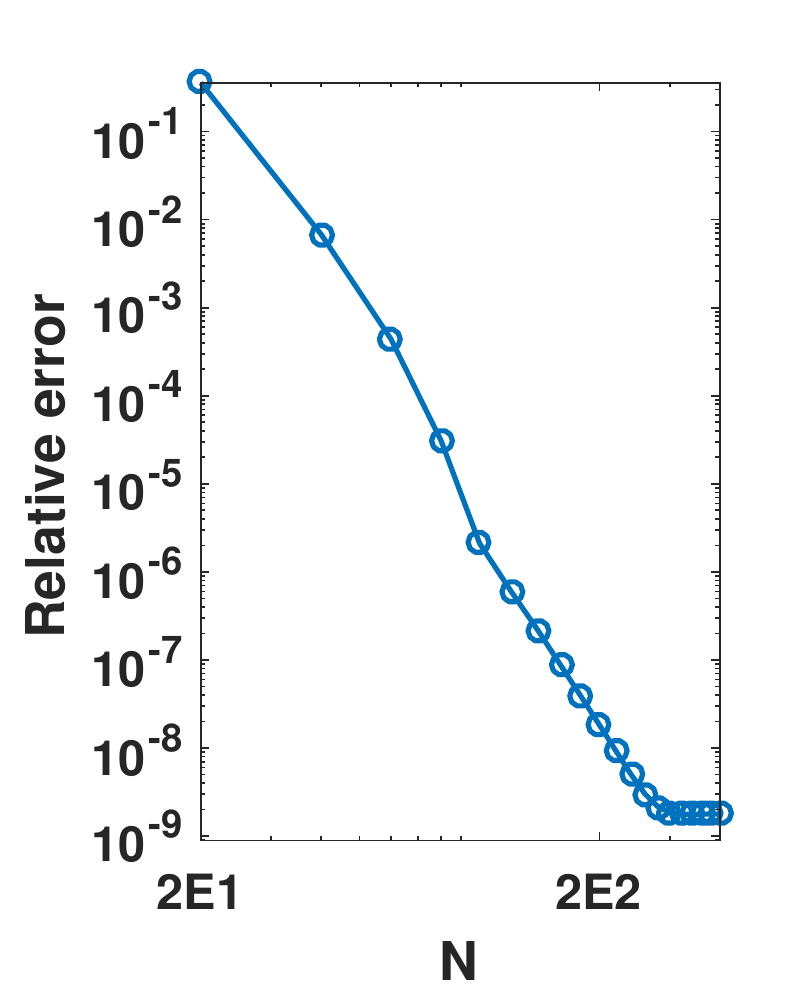}
	(b)\includegraphics[width=0.4\textwidth]{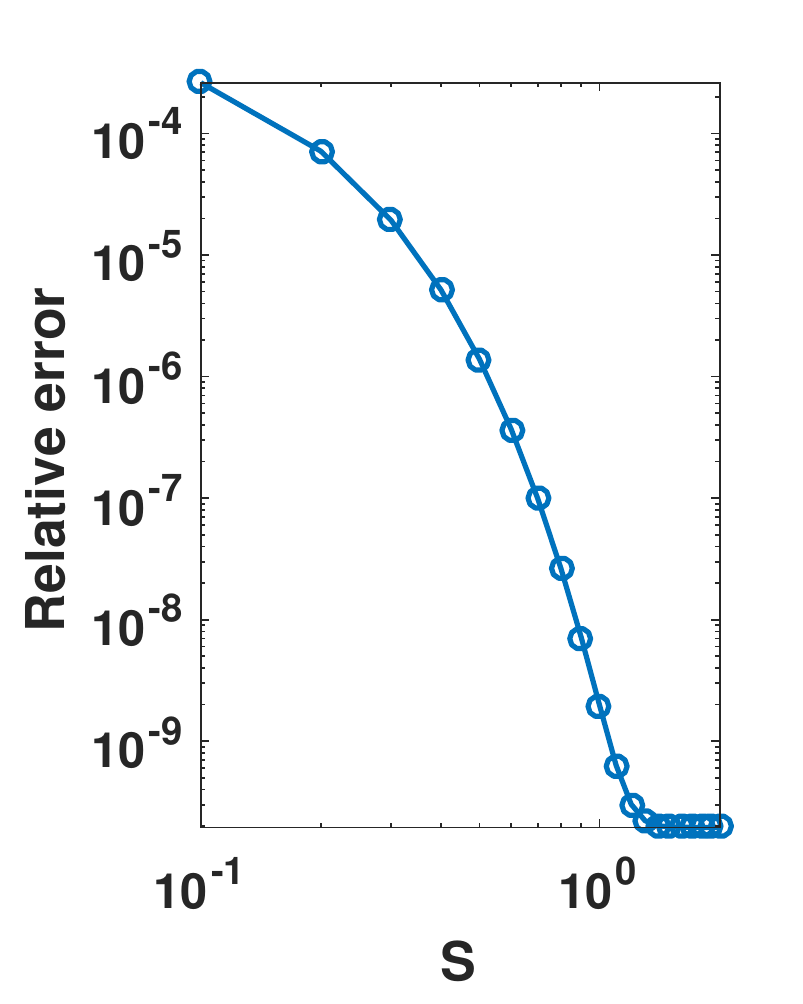}
	\caption{Example 1: Using the exact solution as a reference solution: (a)
		numerical error of $u^{tot}$ on $\Gamma_P$ against total number of
		points $N$ when $S=1$; (b) numerical error of $u^{tot}$ on $\Gamma_P$ against the
		absorbing magnitude $S$ when $N=400$. }
    \label{fig:ex1:cmp3}
\end{figure}

To illustrate the order of accuracy, we study numerical error of $u^{tot}$
on $\Gamma_{P}$ against the number of grid points $N$ in discretizing
$\Gamma_{AB}$ when $S=1$. Since grid points vary for different values of
$N$, we choose to evaluate $u^{tot}$ at grid points on $\Gamma_P$ when
$N=20$, referred to as a reference set of points, to realize the comparison;
for $N\neq 20$, we just interpolate the numerical solution onto the
reference set of points by (\ref{eq:f:trg}). Using the exact solution
$u^{tot}_{exa}$ as a reference solution, we compute numerical errors for
different values of $N$, as depicted in Figure~\ref{fig:ex1:cmp3}(a),
where the vertical axis represents the relative error, the horizontal axis
represents $N$, and both axes are logarithmically scaled. Clearly, slope
of the decaying part of the curve reveals that our method exhibits at least
a seventh-order accuracy.   

To illustrate that our PML effectively terminates the outgoing wave, we now
fix $N=400$ and compute $u^{tot}$ at grid points on $\Gamma_P$ for
different values of $S$, ranging from $0.1$ to $2$; the grid points now are
independent of $S$.  Using the exact solution $u_{exa}^{tot}$ as a
reference solution, we compute relative errors for different values of
$S$, as shown in Figure~\ref{fig:ex1:cmp3}(b), where both axes are
logarithmically scaled.  We observe that the relative error decays
exponentially at the beginning for $S$ in a range of small values, and
however it terminates for larger $S$. We remark that to maintain an
exponentially decaying error for larger $S$, one has to choose larger $N$
to increase the number of points in the PML and to decrease the
discretization error. From Figure~\ref{fig:ex1:cmp3}, we easily see that
the numerical solution for $N=400$ and $S=1$ attains eight significant digits.

To conclude this example, we observe that numerical accuracy in fact can be
improved by two approaches: increasing $N$ and increasing $S$. When exact
solution is not available, it is reasonble that one combines the
convergence curve of relative error against $N$ for a fixed $S$, and the
convergence curve of relative error against $S$ for a fixed $N$ to truly
discover how accurate the solution has obtained, as will be shown below.

\subsection{Example 2: Two semicircles}
\begin{figure}[!ht]
  \centering
  (a)\includegraphics[width=0.4\textwidth]{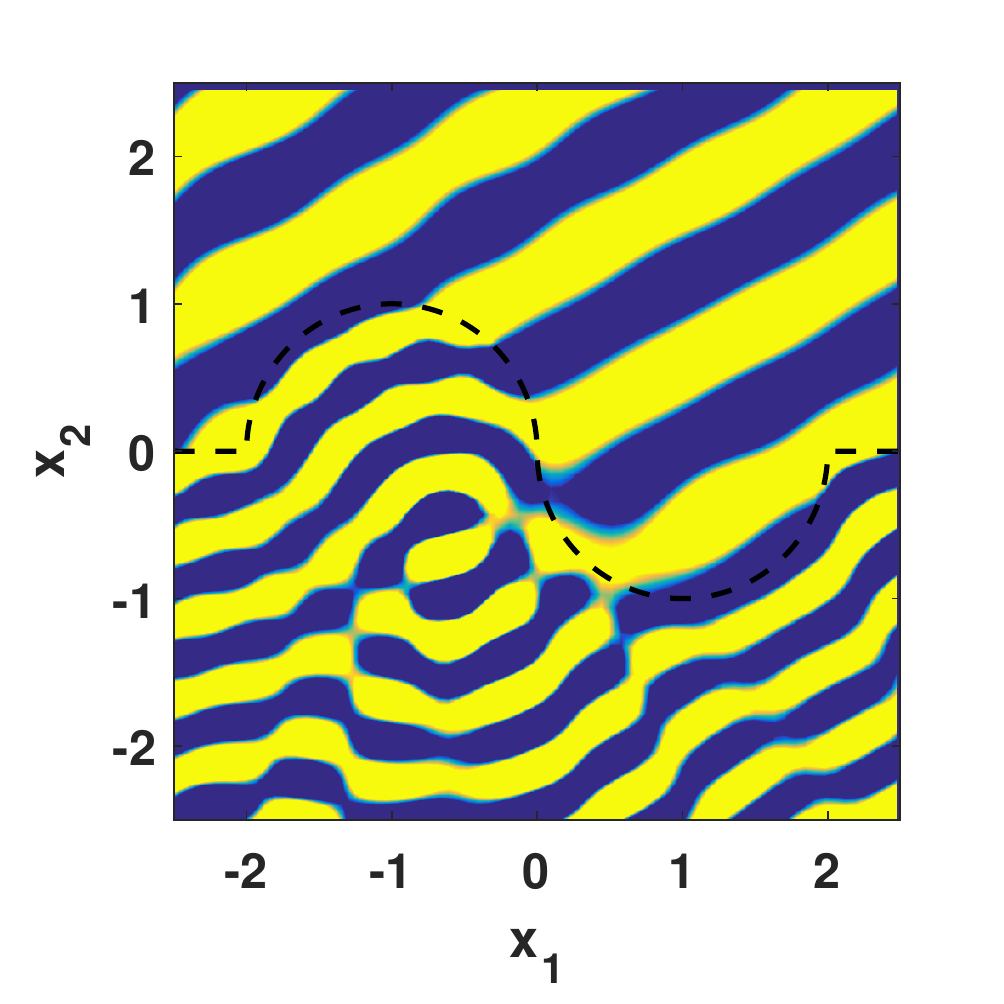}
  (b)\includegraphics[width=0.4\textwidth]{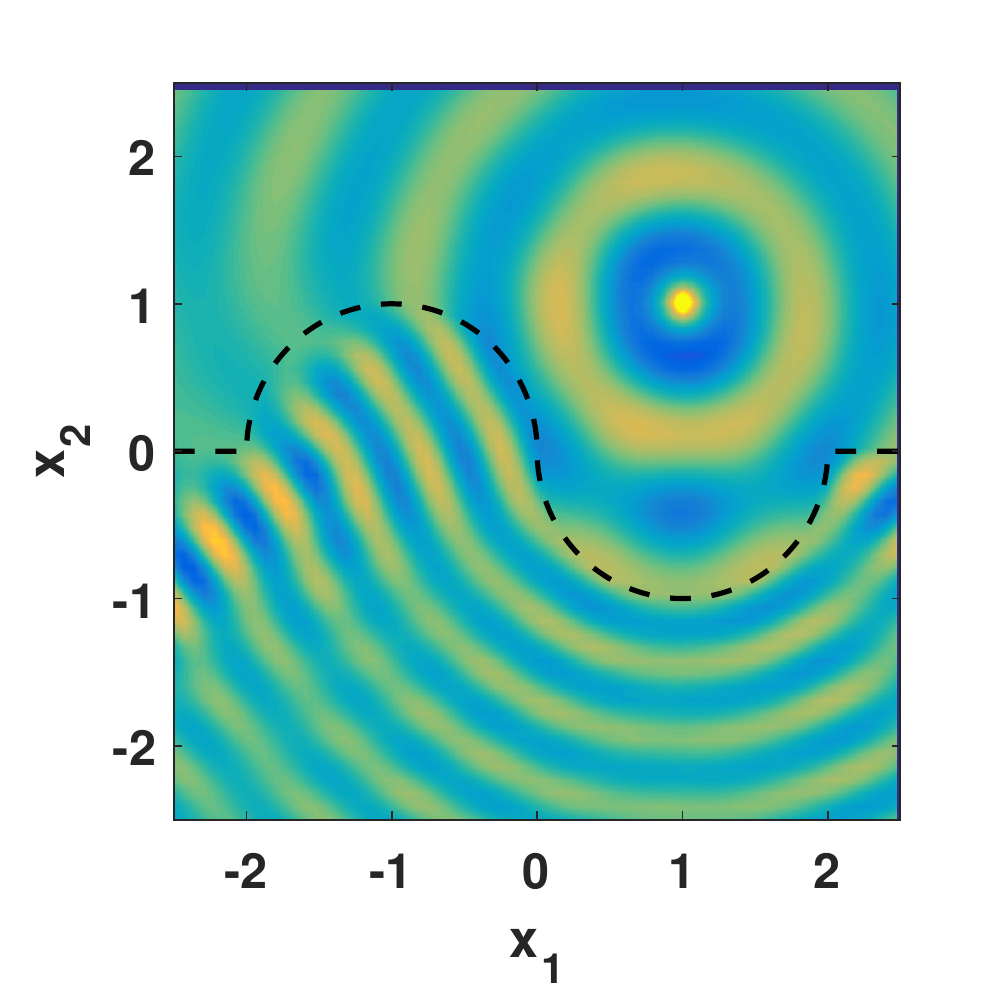}
	\caption{Example 2: real part of $u^{tot}$ on
		$[-2.5,2.5]\times[-2.5,2.5]$. (a) plane incident wave with angle
		$\alpha=\frac{\pi}{3}$; (b) cylindrical wave with source ${\bf
	x}^*=(1,1)$.  Dashed line indicates location of $\Gamma$.   }
    \label{fig:ex2:cmp1}
\end{figure}

In the second example, we consider a local perturbation that consists of two 
connected semicircles of radius $1$; the interface is shown as dotted line
in Figure~\ref{fig:ex2:cmp1}. Suppose again $n_1=1$, $n_2=2$ and $k_0=2\pi$
with wavelength $\lambda=1$. We consider two incident waves: 
\begin{itemize}
	\item[(i)] a plane incident wave with incident angle
		$\alpha=\frac{\pi}{3}$;  
	\item[(ii)] a cylindrical wave due to point source $x^*=(1,1)$.  
\end{itemize}

In the implementation, we take $a_1=2.5$ and $T=1$ so that $\Gamma_P$ becomes $\{(x_1,x_2)|-2.5\leq x_1\leq 2.5\}$ while the PML region is
$\{(x_1,x_2)|2.5\leq|x_1|\leq 3.5\}$. The total wave field $u^{tot}$ for the
two incident waves in $[-2.5,2.5]\times[-2.5,2.5]$ is computed and plotted in
Figure~\ref{fig:ex2:cmp1} (a) and (b), respectively, based on a numerical
solution using $N=1600$ grid points on $\Gamma_{AB}$.  

\begin{figure}[!ht]
  \centering
  (a)\includegraphics[width=0.4\textwidth]{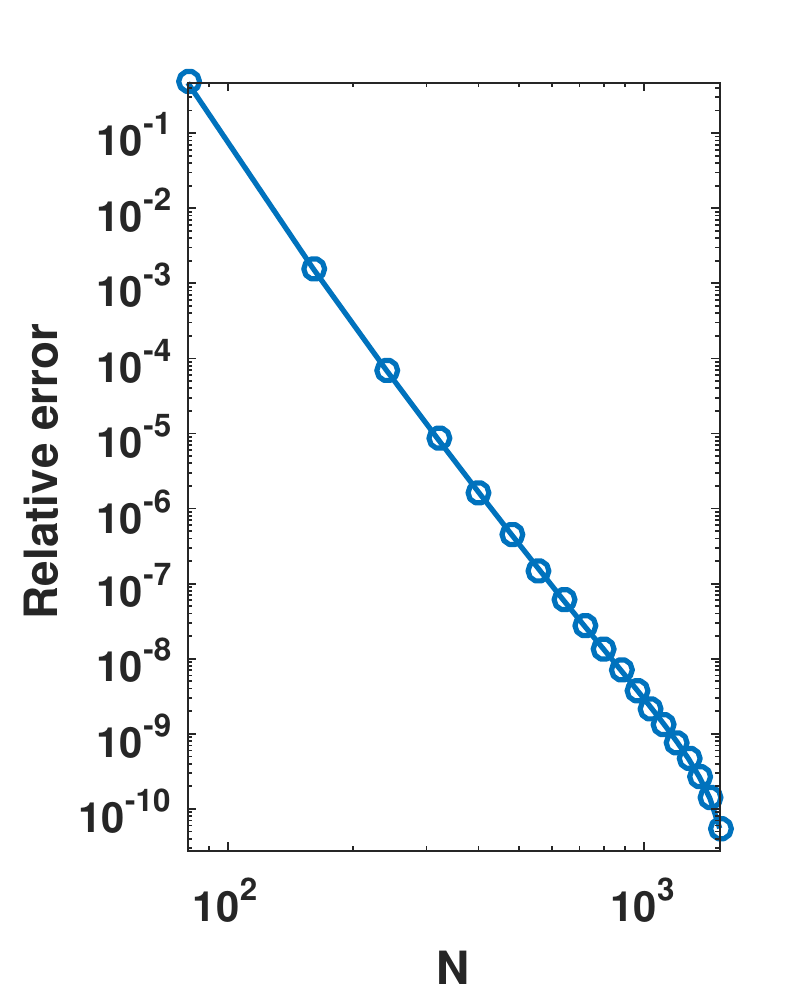}
  (b)\includegraphics[width=0.4\textwidth]{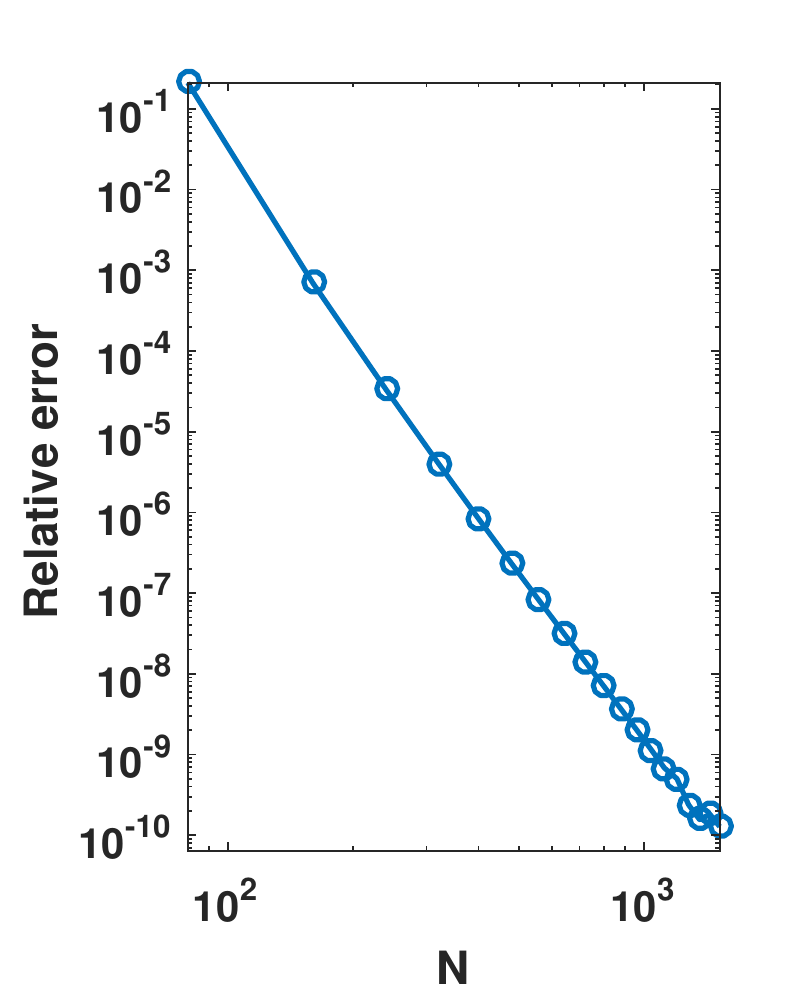}
	\caption{Example 2: numerical error of $u^{tot}$ on $\Gamma_P$ against
		total number of points $N$ when $S=1$: (a) plane incident  wave with angle
		$\alpha=\frac{\pi}{3}$;  (b) cylindrical wave due to point
	source $x^*=(1,1)$.  }
    \label{fig:ex2:cmp2}
\end{figure}

To illustrate the order of accuracy for either incident wave, we compute
numerical error of $u^{tot}$ on $\Gamma_{P}$ against
the number of grid points $N$ when $S=1$. As in example
1, a reference set of points is chosen as the grid points on
$\Gamma_P$ when $N=160$. The reference solution is obtained by computing
$u^{tot}$ at the reference set of points when $N=1600$ grid points are used.
Numerical results for both incident waves are shown in
Figure~\ref{fig:ex2:cmp2},
which shows that our results exhibit a seventh-order accuracy for both
incident waves.  

\begin{figure}[!ht]
  \centering
  (a)\includegraphics[width=0.4\textwidth]{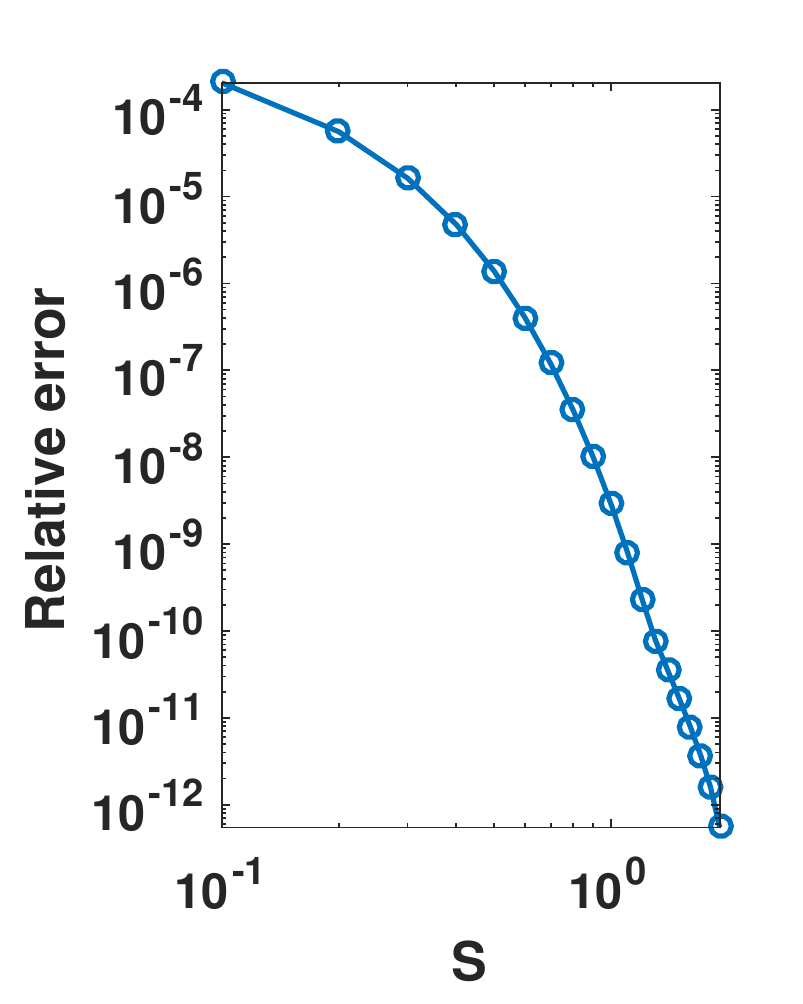}
  (b)\includegraphics[width=0.4\textwidth]{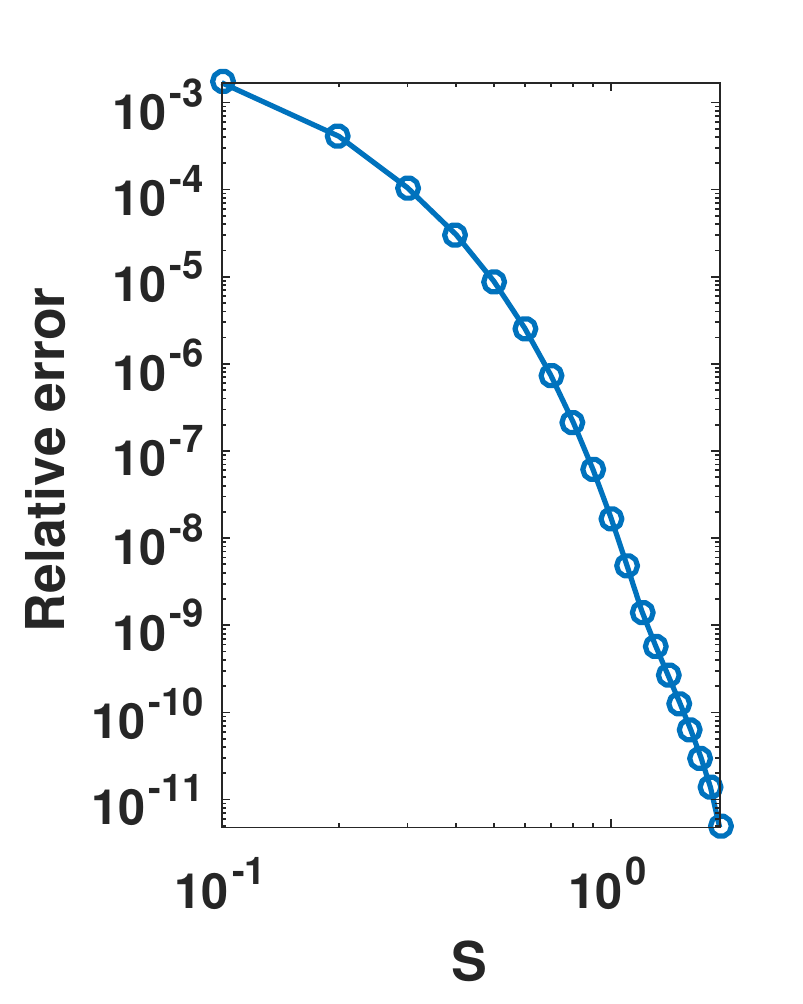}
	\caption{Example 2: numerical error of $u^{tot}$ on $\Gamma_P$ against
		the absorbing magnitude $S$ when $N=1600$: (a) plane incident  wave with angle
		$\alpha=\frac{\pi}{3}$;  (b) cylindrical wave due to point
	source $x^*=(1,1)$.}
    \label{fig:ex2:cmp3}
\end{figure}

To illustrate that our PML effectively terminates the outgoing wave for each
incident wave, we now fix $N=1600$ and compute $u^{tot}$ at grid points on
$\Gamma_P$
 for different values of $S$, ranging
from $0.1$ to $2$; the grid points now are independent of $S$.  Considering
the numerical solution $u^{tot}$ for $S=2$ as a reference solution, we
compute relative errors for different values of $S$. Numerical results are
shown in  Figure~\ref{fig:ex2:cmp3}.
Clearly, we observe that numerical error for each incident wave decays
exponentially at the beginning when $S$ is not very large, and then decays
algebraically for larger $S$ as $N$ is fixed.     

At last, combining Figures~\ref{fig:ex2:cmp2}(a) and \ref{fig:ex2:cmp3}(a),
we see that our numerical solution for the plane incident wave attains
eight significant digits when $N=1600$ and $S=1$.
Similarly, combining Figures~\ref{fig:ex2:cmp2}(b) and
\ref{fig:ex2:cmp3}(b), we see that our numerical solution for the
cylindrical incident wave attains eight significant digits when 
$N=1600$ and $S=1$.   

\subsection{Example 3: An obstacle above the interface}
In this example, we study a more complicated structure, where an obstacle is
placed above the interface. With the obstacle invovled, our PML-based BIE
formulation only requires an extra NtD operator defined on the boundary of the
obstacle, which can be obtained by a regular BIE in physical domain as
described in \cite{lulu14}. Then, according to transmission conditions on
the obstacle and the interface, the final linear system can be obtained with
ease.  

\begin{figure}[!ht]
  \centering
  (a)\includegraphics[width=0.4\textwidth]{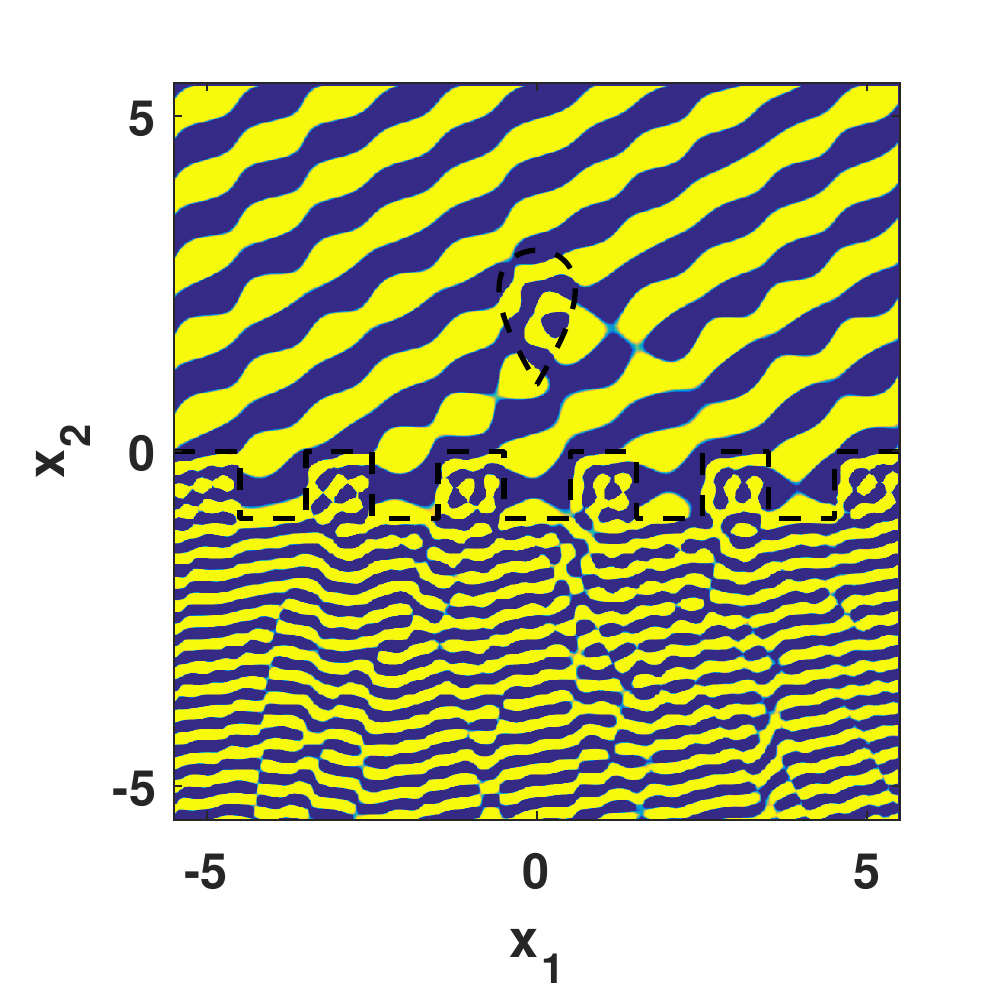}
  (b)\includegraphics[width=0.4\textwidth]{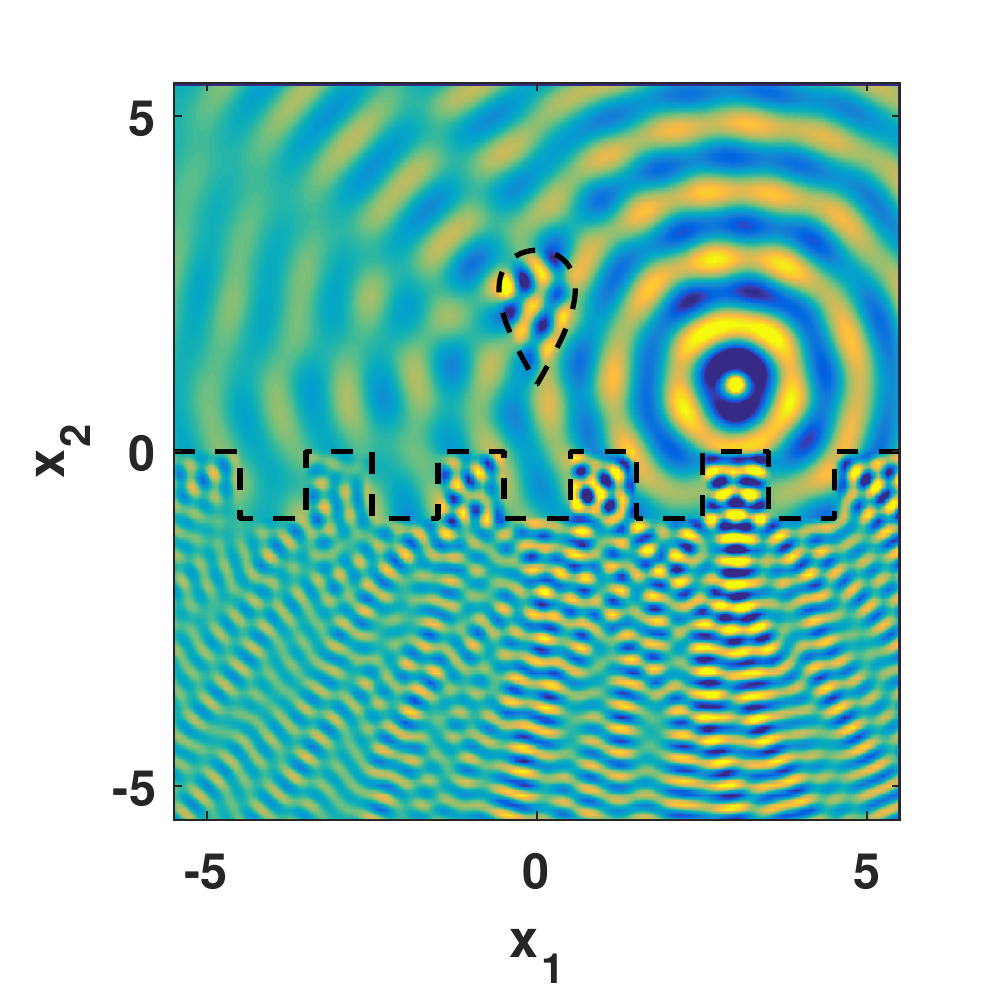}
	\caption{Example 3: real part of $u^{tot}$ on
		$[-5.5,5.5]\times[-5.5,5.5]$: (a) plane incident  wave with angle
		$\alpha=\frac{\pi}{3}$; (b) cylindrical plane wave with point
		source $x^*=(3,1)$. Dashed line indicates location of
		$\Gamma$. }
    \label{fig:ex3:cmp1}
\end{figure}

Suppose refractive index of the obstacle is $n_{ob}=2$, $n_1=1$, $n_2=3$,
and $k_0=2\pi$ with $\lambda = 1$. The basic structure is shown in
Figure~\ref{fig:ex3:cmp1}, where a drop shape is placed one unit above the
interface which contains five uniformly spaced indentations.  We consider two
incident waves: 
\begin{itemize}
	\item[(i)] a plane incident  wave with incident angle
		$\alpha=\frac{\pi}{3}$;  
	\item[(ii)] a cylindrical wave due to point source $x^*=(3,1)$.  
\end{itemize}

In the implementation, we take $a_1 = 5.5$ and $T = 1$ so that $\Gamma_P$ becomes $\{(x_1,x_2)|-5.5\leq x_1\leq 5.5\}$, while the PML
domain becomes $\{(x_1,x_2)|5.5\leq|x_1|\leq 6.5\}$. The total wave field
$u^{tot}$ for the two incident waves in $[-5.5,5.5]\times[-5.5,5.5]$ is
computed and plotted in Figure~\ref{fig:ex3:cmp1} (a) and (b), respectively,
based on a numerical solution on $\Gamma_{AB}$ and the obstacle boundary
$\Gamma_{ob}$, using $N=3150$ grid points on $\Gamma_{AB}$ ($150$ points per
segment) and $N_{ob}=800$ grid points on $\Gamma_{ob}$, the boundary of the
obstacle.  

\begin{figure}[!ht]
  \centering
  (a)\includegraphics[width=0.4\textwidth]{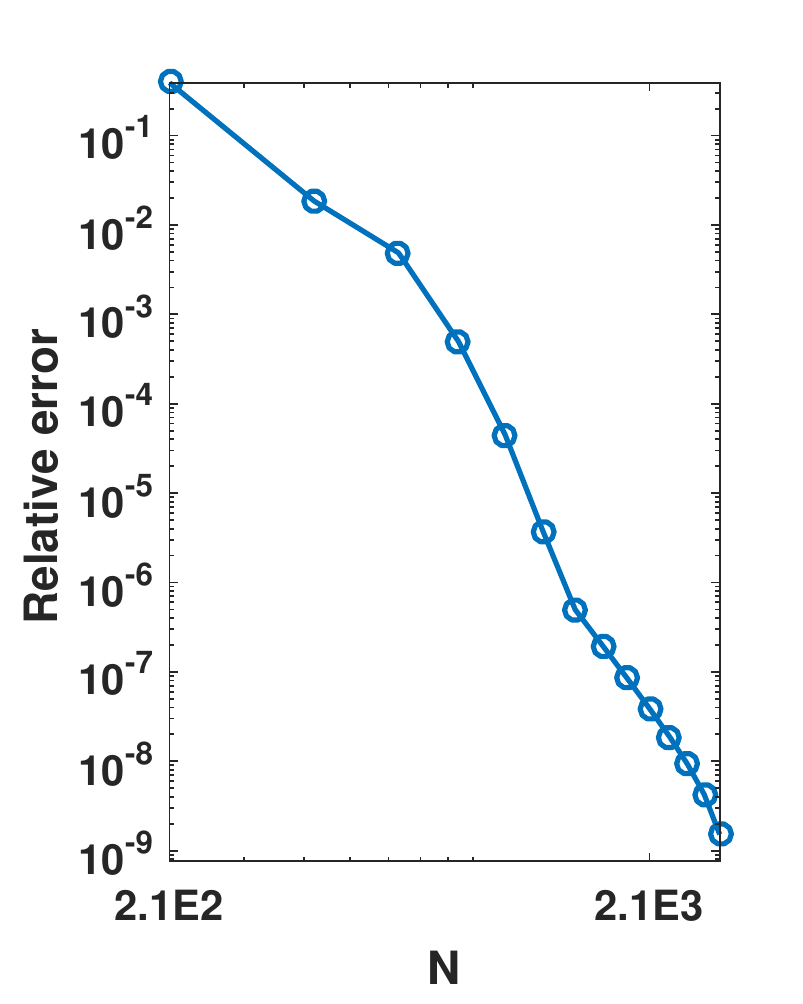}
  (b)\includegraphics[width=0.4\textwidth]{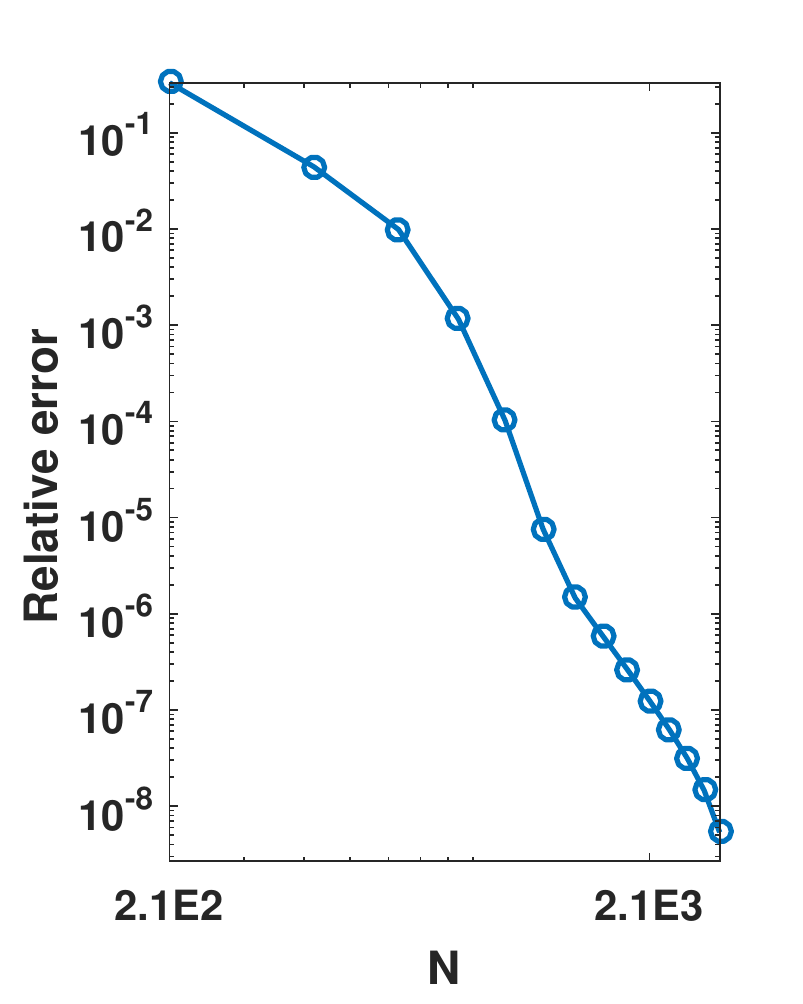}
	\caption{Example 3: numerical error of $u^{tot}$ on $\Gamma_P$ against
		total number of points $N$ when $S=1$: (a) plane incident  wave with angle
		$\alpha=\frac{\pi}{3}$;  (b) cylindrical wave due to a point
		source $x^*=(3,1)$. }
    \label{fig:ex3:cmp2}
\end{figure}

To illustrate the order of accuracy for either incident wave, we study
numerical error of $u^{tot}$ on $\Gamma_{P}$ against
the number of grid points $N$ in discretizing $\Gamma_{AB}$ when $S=1$, where we fix
the number of grid points on $\Gamma_{ob}$ to be $N_{ob}=800$. As in example
1, a reference set of points is chosen as the grid points on $\Gamma_P$ when
$N=840$. The reference solution is obtained by computing $u^{tot}$ at the
reference set of points when $N=3150$ grid points are used. Numerical
results for both incident waves are shown in Figure~\ref{fig:ex3:cmp2},
which shows that our results roughly exhibit a seventh-order accuracy for
both incident waves.  

\begin{figure}[!ht]
  \centering
  (a)\includegraphics[width=0.4\textwidth]{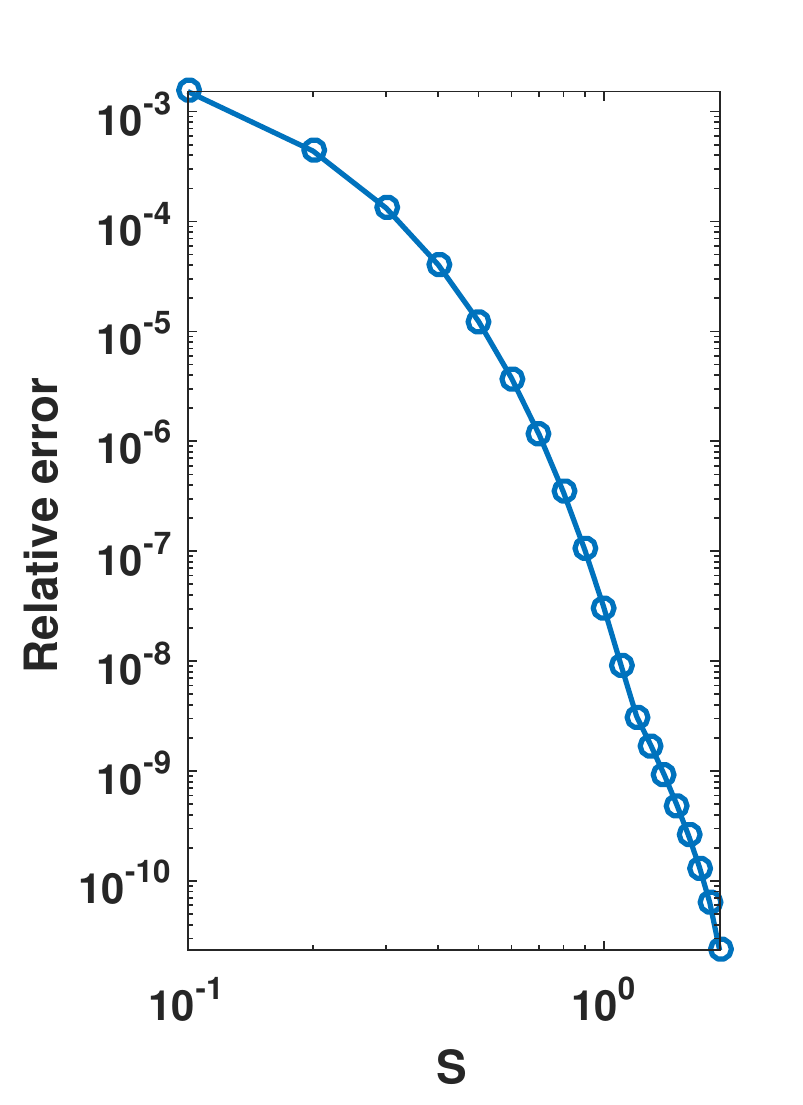}
  (b)\includegraphics[width=0.4\textwidth]{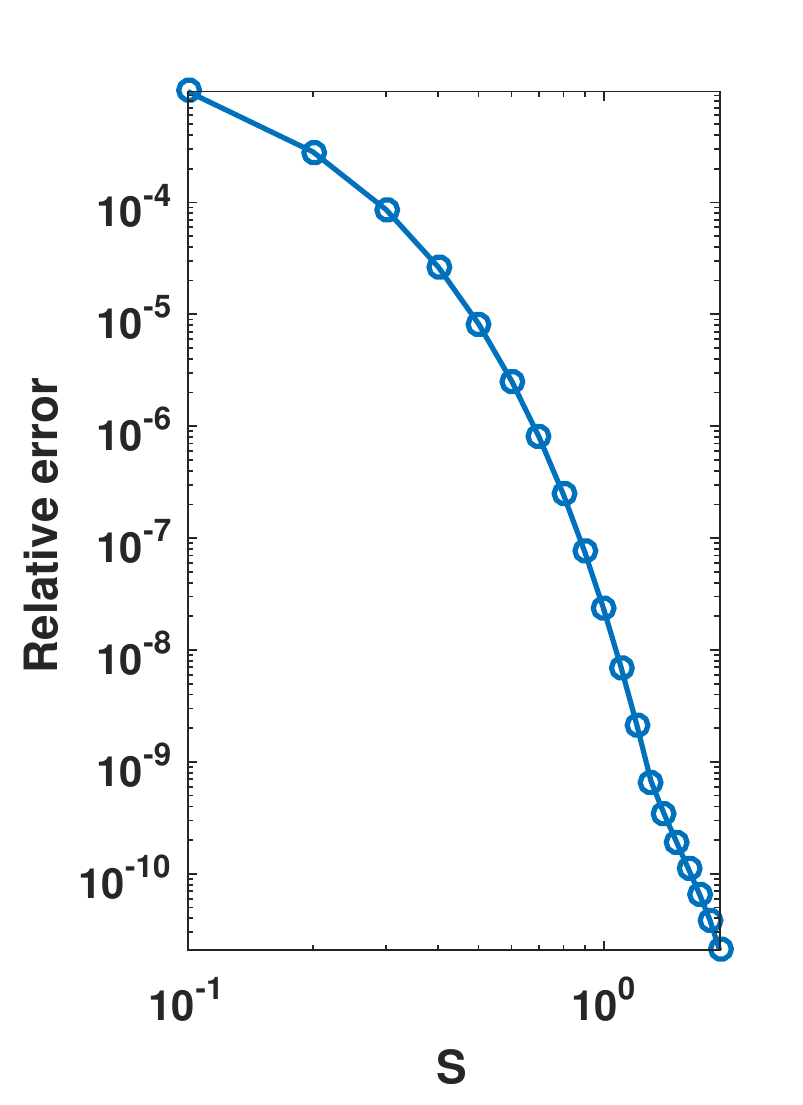}
	\caption{Example 3: numerical error of $u^{tot}$ on $\Gamma_P$ against
		the absorbing magnitude $S$ when $N=3150$: (a) plane incident  wave with angle
		$\alpha=\frac{\pi}{3}$;  (b) cylindrical wave due to a point
		source $x^*=(3,1)$. }
    \label{fig:ex3:cmp3}
\end{figure}

To illustrate that our PML effectively terminates the outgoing wave for each
incident wave, we now fix $N=3150$ and compute $u^{tot}$ at grid points on
$\Gamma_{P}$ for different values of $S$, ranging
from $0.1$ to $2$; the grid points now are independent of $S$.  Considering
the numerical solution $u^{tot}$ for $S=2$ as a reference solution, we
compute relative errors for different values of $S$. Numerical results are
shown in  Figure~\ref{fig:ex3:cmp3}.
Clearly, we observe that numerical error for each incident wave decays
exponentially at the beginning when $S$ is not very large, and then decays
algebraically for larger $S$ as $N$ is fixed.     

At last, combining Figures~\ref{fig:ex3:cmp2}(a) and \ref{fig:ex3:cmp3}(a),
we see that our numerical solution for the plane incident wave attains 
seven significant digits when $N=3150$ and $S=1$. Similarly, combining
Figures~\ref{fig:ex3:cmp2}(b) and \ref{fig:ex3:cmp3}(b), we see that our
numerical solution for the cylindrical incident wave attains seven
significant digits when $N=3150$ and $S=1$.

\subsection{Example 4: Interface with different elevations at infinity}

In previous examples, flat part of the interface away from the local
perturbation $P$ have the same elevations at infinity. However, if the flat part has
different elevations toward infinity, then all existing methods based on
layered Green's function break down since now for the background layered
medium, an explicit form of the layered medium Green's function in terms of
Sommefeld integrals is hard to develop. To conclude this section, we study
such a challenging example.   

For a plane incident wave, using a flat part on one side (left or right)  to define $u^{tot}_0$
can only suppress the reflective and transmittive waves in $u^s$ on the same
side but not on the other side, since the reflection and transmission
coefficients are different on each side. Consequently, it is possible that
the difference field $u^{tot}-u^{tot}_0$ is not outgoing in all directions,
e.g., if $u^{inc}$ a normal incident wave. The current PML-based BIE
formulation fails in this case. We expect to address this issue in an
ongoing project.   

\begin{figure}[!ht]
  \centering
  (a)\includegraphics[width=0.4\textwidth]{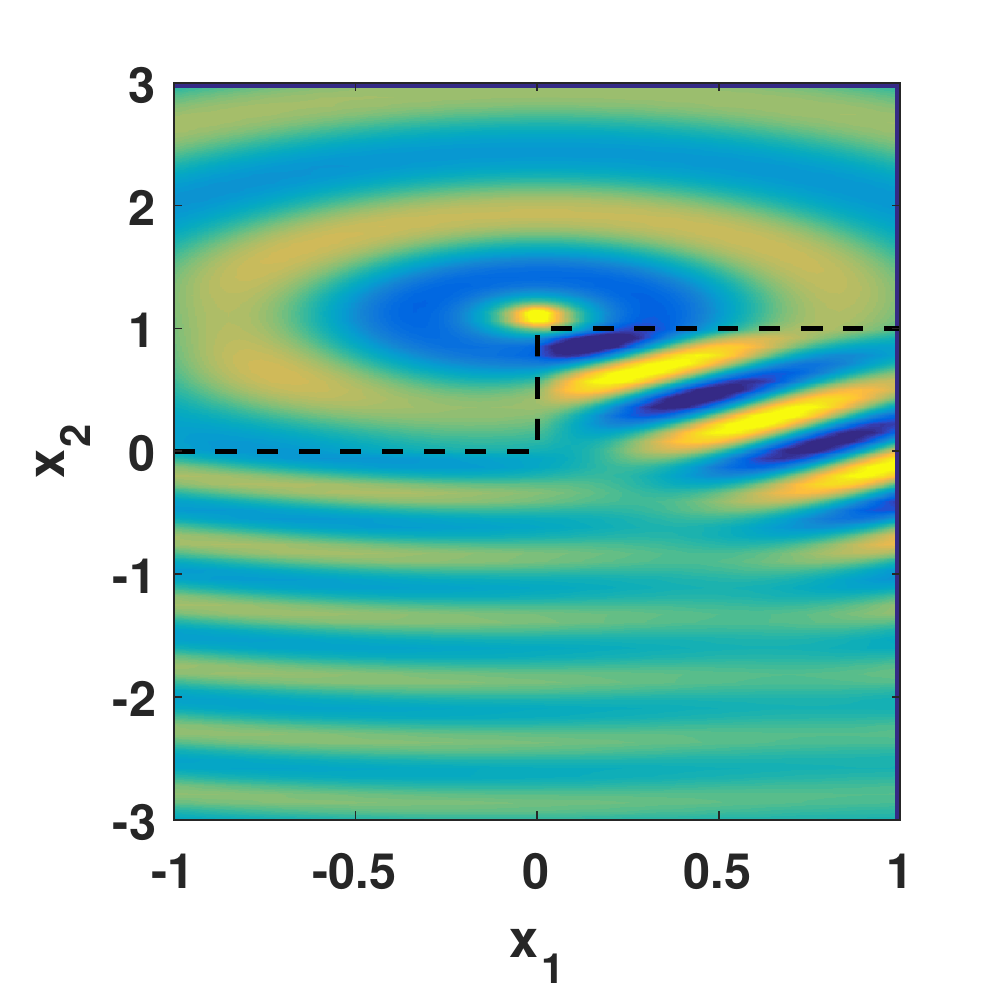}
	\caption{Example 4: real part of $u^{tot}$ on
		$[-1,1]\times[-3,3]$ for cylindrical plane wave due to
		source ${\bf r}_0=[0,1.1]^{T}$. Dashed line indicates location of
		$\Gamma$. }
    \label{fig:ex4:cmp1}
\end{figure}

Fortunately, when $u^{inc}$ is a cylindrical wave due to a point source, we
may still use $u^{tot}_0$ defined in (\ref{eq:utot0:2}) to construct an
outgoing wave $u^s$ such that our PML-based BIE formulation still works.
To justify the methodology, we test a very simple structure where two
half-lines with different elevations are connected just by a line segment of
$1$ unit, as shown in Figure~\ref{fig:ex4:cmp1}, where we suppose $n_1=1$,
$n_2=2$ and $k_0=2\pi$ with wavelength $\lambda=1$. We consider a
cylindrical incident wave due to point source $x^*=(0,1.1)$.   

In the implementation, we take $a_1 = 1$ and $T = 1$ so that $\Gamma_P$
becomes $\{(x_1,x_2)|-1\leq x_1\leq 1\}$, while the PML domain becomes
$\{(x_1,x_2)|1\leq|x_1|\leq 2\}$. The total wave field $u^{tot}$ for the
incident wave in $[-1,1]\times[3,3]$ is computed and plotted in
Figure~\ref{fig:ex4:cmp1},
based on a numerical solution using $N=2400$ grid points on $\Gamma_{AB}$
($800$ points per smooth segment).  

\begin{figure}[!ht]
  \centering
  (a)\includegraphics[width=0.4\textwidth]{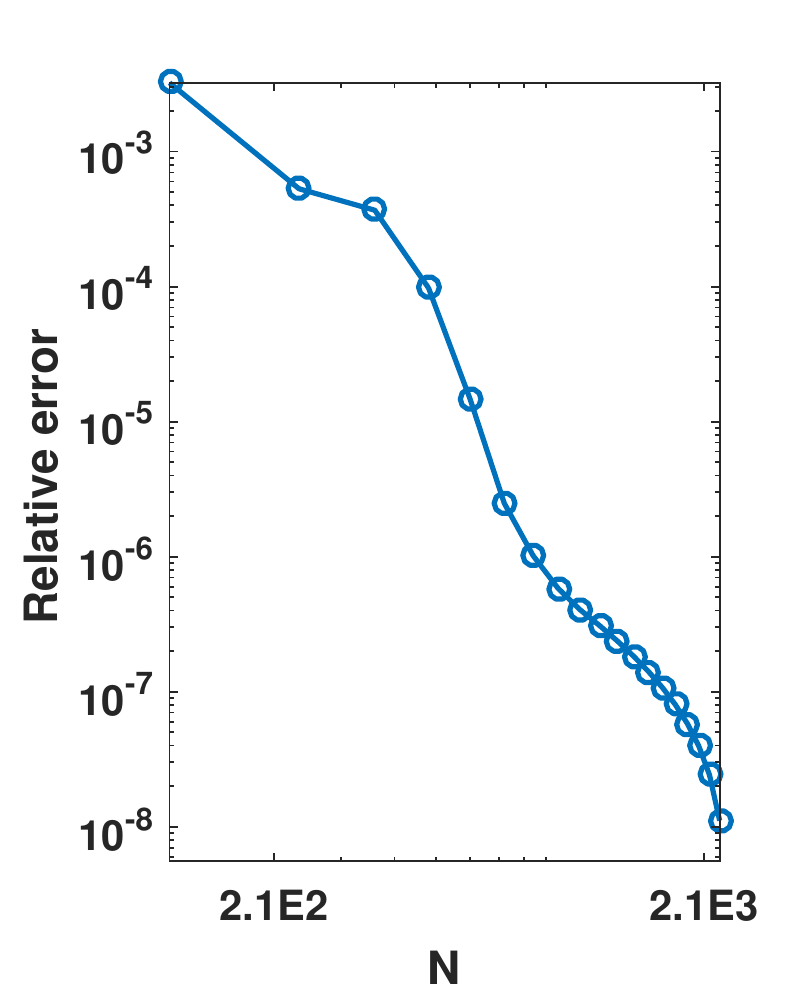}
  (b)\includegraphics[width=0.4\textwidth]{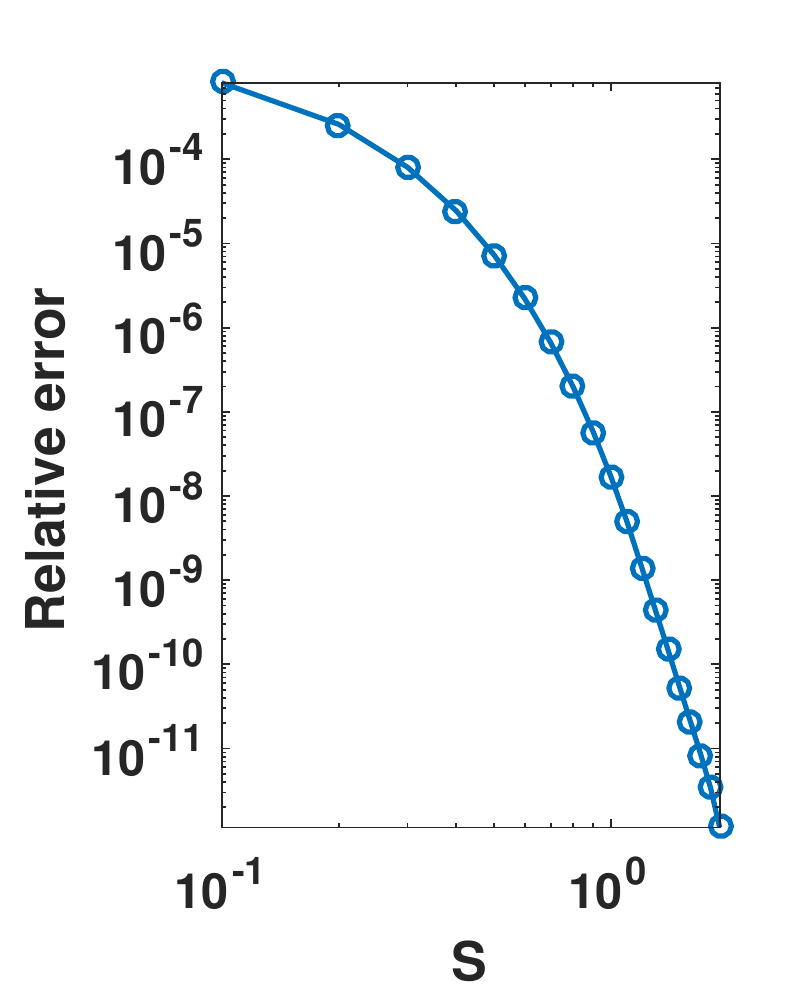}
	\caption{Example 4: (a) numerical error of $u^{tot}$ on $\Gamma_P$ against
		the total number of points $N$ when $S=1$; (b) numerical error of $u^{tot}$
		on $\Gamma_P$ against the absorbing magnitude $S$ when $N=2400$.  }
    \label{fig:ex4:cmp2}
\end{figure}

To illustrate the order of accuracy, we study
numerical error of $u^{tot}$ on $\Gamma_{P}$ against
the number of grid points $N$ in discretizing $\Gamma_{AB}$ when $S=1$. As in example
1, a reference set of points is chosen as the grid points on $\Gamma_{P}$
when $N=120$. The reference solution is
obtained by computing $u^{tot}$ at the reference set of points when $N=2400$ grid
points are used. Numerical results are shown in
Figure~\ref{fig:ex4:cmp2}(a),
which shows that our results roughly exhibit a fourth-order accuracy.  

To illustrate that our PML effectively terminates the outgoing wave, we now
fix $N=2400$ and compute $u^{tot}$ at grid points on $\Gamma_{P}$ for
different values of $S$, ranging from $0.1$ to $2$; the grid points now are
independent of $S$.  Considering the numerical solution $u^{tot}$ for $S=2$
as a reference solution, we compute relative errors for different values of
$S$. Numerical results are shown in Figure~\ref{fig:ex4:cmp2} (b).  Clearly,
we observe that numerical error decays exponentially at the beginning when
$S$ is not very large, and then decays algebraically for larger $S$ as $N$
is fixed.

Finally, combining Figures~\ref{fig:ex4:cmp2}(a) and \ref{fig:ex4:cmp2}(b),
we see that our numerical solution for the cylindrical wave attains seven
significant digits when $N=2400$ and $S=1$. 

\section{Conclusion}
For 2D scattering problems in layered media with unbounded interfaces, we
developed a PML-based BIE method that relies on the Green's function of
PML-transformed free space.  The method avoid the difficulty of evaluating
the expensive Sommerfeld integrals in common BIE methods based on Green's
functions of layered media. Similar to other BIE methods based on the free
space Green's function, integral equations are formulated on unbounded
interfaces of the background media and these interfaces  must be truncated.
Although existing methods such as the windowing function method
\cite{brudel14,mon08,brulyoperaratur16,laigreone15}, are also effective in
truncating interfaces, our method is particularly simple, since the
truncation simply follows the well-established PML technique.  Notice that
the Green's function of PML-transformed free space is simply obtained from
the usual Green's function by extending the argument to complex space, and
it is very easy to evaluate. 

Since our main purpose is to develop a PML-based method and demonstarte its
effectiveness for truncating the unbounded interfaces, we have used a
simple BIE formulation involving  the single- and double-layer boundary
integral operators only. In addition, we used the DtN maps to simplify the
final linear system. Numerical examples are presented for scattering
problems involving two homogeneous media separated by an interface with
local perturbations, and possibly with additional obstacles.  The integral
equations are discretized using  a graded mesh technique,  Alpert's  sixth
order hybrid Gauss-trapezoidal rule for logarithmic singularities, and a
stabilizing technique.  Numerical results indicate that the truncation of
interfaces by PML is highly effective, and  accurate solutions can be
obtained using PMLs with a thickness of one wavelength. 

Although our current implementation is somewhat limited,  the PML-based BIE
method can be extended in a number of directions.  Obviously, the method
can be used to study scattering problems in multi-layered media  with local
perturbations, embedded obstacles and penetrable structures.  Besides
scattering problems, the method can also be used to study eigenvalue
problems, such as the problem for guided modes in open waveguide
structures.  We are planning to address some of these problems in our
future works.  

\section*{Acknowledgement}
Y. Y. Lu is partially supported by the Research Grants Council of Hong
Kong Special Administrative Region, China (Grant No. CityU 11301914). J.
Qian is partially supported by NSF grants 1522249 and 1614566. 
\begin{appendix}
\section*{Appendix}

In this appendix, we will show that equation (\ref{eq:K_0:1}) holds for any 
$x^o=[x_1^o,x_2^o]$ on $\Gamma_{AB}$.   

At first, using the Green's representation theorem, we easily see that
\begin{equation}
	\label{eq:K_0:limit}
	\tilde{\cal K}_0[1](x^o) = \lim_{\varepsilon\rightarrow
	0^+}2\int_{\partial B(x^o,\varepsilon)\cap \bar{\Omega}} \partial_{{\bm \nu}_c}
	G_0(x^o,x) 1 ds(x),
\end{equation}
where $\partial B(x^o,\varepsilon)$ is the boundary of circle
$B(x^o,\varepsilon)$ of radius $\varepsilon$ centered at $x^o$, and the unit
normal vector ${\bm \nu}$ now points toward $\Omega$.    

Thus for sufficiently small $\varepsilon$, one can parameterize
$\partial B(x^o,\varepsilon)\cap\bar{\Omega}$ by $x = x^o+\varepsilon(\cos t,\sin
t)$ for $t\in[\theta_1,\theta_2]$ where the inner angle
$\theta=\theta_2-\theta_1$.   
 
Clearly, according to its definition (\ref{eq:def:K0}), we can discretize
$\tilde{\cal K}_0$ as
\begin{equation}
	\label{eq:dis:K0}
	\tilde{\cal K}_0[1](x^o) = -\frac{1}{\pi}\lim_{\varepsilon\rightarrow 0+}\int_{\theta_1}^{\theta_2}
	\frac{(\tilde{x}_1 -
	\tilde{x}_1^o)\tilde{x}_2'-\tilde{x}_1'(\tilde{x}_2-\tilde{x}_2^o)}{|\tilde{x}^o-\tilde{x}|^2}
	d t.  
\end{equation}

By definitions of complex stretched coordinates transformation
(\ref{eq:x}), on the boundary $\partial B(x^o,\varepsilon)\cap\bar{\Omega}$,
we have
\begin{align}
  \tilde{x}_1 - \tilde{x}_1^o &= \int_{x_1^o}^{x_1} \alpha_1(s)
  ds\nonumber\\ 
  &= \int_{x_1^o}^{x_1^o+\varepsilon \cos t}
  \alpha_1(s)ds\nonumber\\ 
  &= \alpha_1(x_1^o)\varepsilon\cos t + O(\varepsilon^2),
\end{align}
and similarly,
\begin{align}
\tilde{x}_2 - \tilde{x}_2^o = \alpha_2(x_2^o)\varepsilon\sin t + O(\varepsilon^2).
\end{align}
Thus,
\begin{align}
	\label{eq:dis:K0:2}
        \tilde{\cal K}_0[1](x^o)&=-\frac{1}{\pi}\lim_{\varepsilon\rightarrow
        0+}\int_{\theta_1}^{\theta_2}\frac{ \alpha_1(x_1^o)
        \alpha_2(x_2^o)\varepsilon^2 +
      O(\varepsilon^3)}{\alpha_1^2(x_1^o)\varepsilon^2\cos^{2}t +
      \alpha_2^2(x_2^o)\varepsilon^2\sin^2 t + O(\varepsilon^3)}
      dt\nonumber\\
&=-\frac{1}{\pi}\int_{\theta_1}^{\theta_2}
        \frac{\alpha_1(x_1^o)\alpha_2(x_2^o)}{\alpha_1^2(x_1^o)\cos^2 t +
        \alpha_2^2(x_2^o)\sin^2 t} dt.
\end{align}
Clearly, if $x^o$ is outside the PML so that $\alpha_1(x^o)=\alpha_2(x^o)=1$, then 
\[
\tilde{\cal K}_0[1](x^o)=-\frac{\theta_2-\theta_1}{\pi}=-\frac{\theta}{\pi}.
\]
If $x^o$ is inside the PML so that $x^o$ is just a smooth point away from the
perturbation curve $P$, then we easily see that $\theta_1=0$, 
$\theta_2=\pi$, and the inner angle $\theta= \pi$. In this case,
\begin{align}
  \label{eq:dis:K0:3}
  \tilde{\cal K}_0[1](x^o)
  &=-\frac{1}{\pi}\int_{0}^{\pi}\frac{\alpha_2(x_2^o)/\alpha_1(x_1^o)\sec^2
  t}{\left(\alpha_2(x_2^o)/\alpha_1(x_1^o)\tan t\right)^2 + 1} dt\nonumber\\
  &=-\frac{1}{\pi}\left(\int_{0}^{\pi/2} d(\arctan \left(\alpha_2/\alpha_1
\tan t\right)) + \int_{\pi/2}^{\pi} d(\arctan \left(\alpha_2/\alpha_1
\tan t\right))\right) \nonumber\\
&=
-\frac{1}{\pi}\left(\frac{\pi}{2}-0\right)-\frac{1}{\pi}\left(0-(-\frac{\pi}{2})\right)
= -\frac{\pi}{\pi}.
\end{align}

\end{appendix}
\bibliographystyle{plain}
\bibliography{wt}
\end{document}